\documentclass[12pt,reqno]{amsart}

\usepackage{latexsym}
\usepackage{amsmath}
\usepackage{amscd}
\usepackage{amssymb}
\usepackage{mathrsfs}
\usepackage{mathtools}
\usepackage{comment}
\usepackage{color,framed}
\usepackage{enumerate}
\usepackage{soul}
\usepackage{stmaryrd}
\usepackage[ruled]{algorithm2e} 
\usepackage{braket}

\usepackage{graphicx}
\usepackage{subfigure}
\usepackage[colorlinks]{hyperref}
\usepackage{url}

\def\R{{\mathbb R}}

\def\Exp{{\rm Exp}}
\def\Log{{\rm Log}}
\def\GL{{\rm GL}}

\DeclareMathOperator{\argmin}{argmin}

\newcommand{\Symp}{\mathrm{Sym}^+}
\newcommand{\vol}{\mathrm{vol}}

\newtheorem{theorem}{Theorem}[section]

\newtheorem{remark}[theorem]{Remark}

\numberwithin{equation}{section}

\begin{document}

\title[An Infinitesimal Probabilistic Model for Manifold PCA]{An Infinitesimal Probabilistic Model for Principal Component Analysis of Manifold Valued Data}


\author{Stefan Sommer}


\address{
Department of Computer Science (DIKU), University of Copenhagen,
  DK-2100 Copenhagen E, Denmark
}


\maketitle

\begin{abstract}
  We provide a probabilistic and infinitesimal view of how the principal component analysis procedure (PCA) can be generalized to analysis of nonlinear manifold valued data. Starting with the probabilistic PCA interpretation of the Euclidean PCA procedure, we show how PCA can be generalized to manifolds in an intrinsic way that does not resort to linearization of the data space. The underlying probability model is constructed by mapping a Euclidean stochastic process to the manifold using stochastic development of Euclidean semimartingales. The construction uses a connection and bundles of covariant tensors to allow global transport of principal eigenvectors, and the model is thereby an example of how principal fiber bundles can be used to handle the lack of global coordinate system and orientations that characterizes manifold valued statistics. We show how curvature implies non-integrability of the equivalent of Euclidean principal subspaces, and how the stochastic flows provide an alternative to explicit construction of such subspaces. We describe estimation procedures for inference of parameters and prediction of principal components, and we give examples of properties of the model on embedded surfaces.

\keywords{principal component analysis, manifold valued statistics, stochastic development, probabilistic PCA, anisotropic normal distributions, frame bundle}
\end{abstract}

\section{Introduction}
A central problem in the formulation of statistical methods for analysis of data in nonlinear spaces is the lack of global coordinate systems and global orientation fields. As an example, consider generalizing the notion of covariance matrix to manifold valued random variables: While the Euclidean definition takes the expectation $E[x^i-E[x]^i]E[x^j-E[x]^j]$ of the product of the coordinate components $x^i$ of the centered random variable $x-E[x]$, the coordinate components are not meaningful in the nonlinear situation as the coordinates themselves are not defined. This fact fundamentally questions what constitutes a natural generalization of covariance. As a second example, consider a standard Euclidean linear latent variable model
\begin{equation}
  y=m+Wx+\epsilon
  \label{eq:linreg}
\end{equation}
on $\R^d$ with mean $m$, coefficient matrix $W$, latent variables $x$, and noise $\epsilon$. The columns of $W$ can be seen as encoding the direction in the Euclidean space connected to a change of each element of $x$. However, on a manifold $M$, $W$ has a priori only meaning for infinitesimal changes $\partial_x$ in the tangent space $T_mM$, and the lack of global orientation prevents a direct translation between such infinitesimal changes, finite perturbations of $m$, and global directions on $M$.

The aim of this paper is to construct a nonlinear manifold generalization of the inherently linear principal component analysis (PCA) procedure, a generalization that is intrinsically based on the geometry of the manifold $M$ and does not resort to a linear approximation of the geometry. The model is based on the Euclidean probabilistic principal component analysis procedure (PPCA, \cite{tipping_probabilistic_1999}) that interprets PCA as a latent variable model \eqref{eq:linreg} with $W$ having low rank $k\le d$. We use the PPCA approach with a probability model based on a notion of infinitesimal covariance and thereby avoid linearizing the nonlinear data space while intrinsically incorporating the effect of data anisotropy, here difference in the principal eigenvalues. The model is related to the probabilistic principal geodesic analysis (PPGA, \cite{zhang_probabilistic_2013}) procedure, however using the probability model and normal distributions defined in \cite{sommer_anisotropic_2015,sommer_modelling_2017}. This construction in particular emphasizes the role of the connection on the manifold in linking infinitesimally close tangent spaces. We refer to the method as being \emph{infinitesimal probabilistic} because the connection allows sequences of random, infinitesimal steps to generate the data probability model.

As a second aim, we wish to exemplify how the use of fiber bundle structures provides a way around the lack of coordinates and global orientations on $M$. The construction in \cite{sommer_anisotropic_2015,sommer_modelling_2017} essentially enlarges the manifold by equipping it with a structure group at each point and hence a principal fiber bundle structure. An example of this is the frame bundle $FM$, viewed as the bundle of invertible linear maps $\GL(\R^d,TM)$, but we will also encounter lower-rank versions $F^kM$ of $FM$, and the quotient bundle $\Symp$ of symmetric positive tensors on $TM$. Elements of these bundles are here used to model the local anisotropy and covariance of the data. The incorporation of nontrivial covariance couples with curvature leading to families of paths that extend geodesics as being, in a certain sense, most probable paths between data points \cite{sommer_anisotropically_2016}. 

\begin{figure}[ht]
    \centering
    \includegraphics[width=.48\columnwidth, trim = 130 150 200 150,clip]{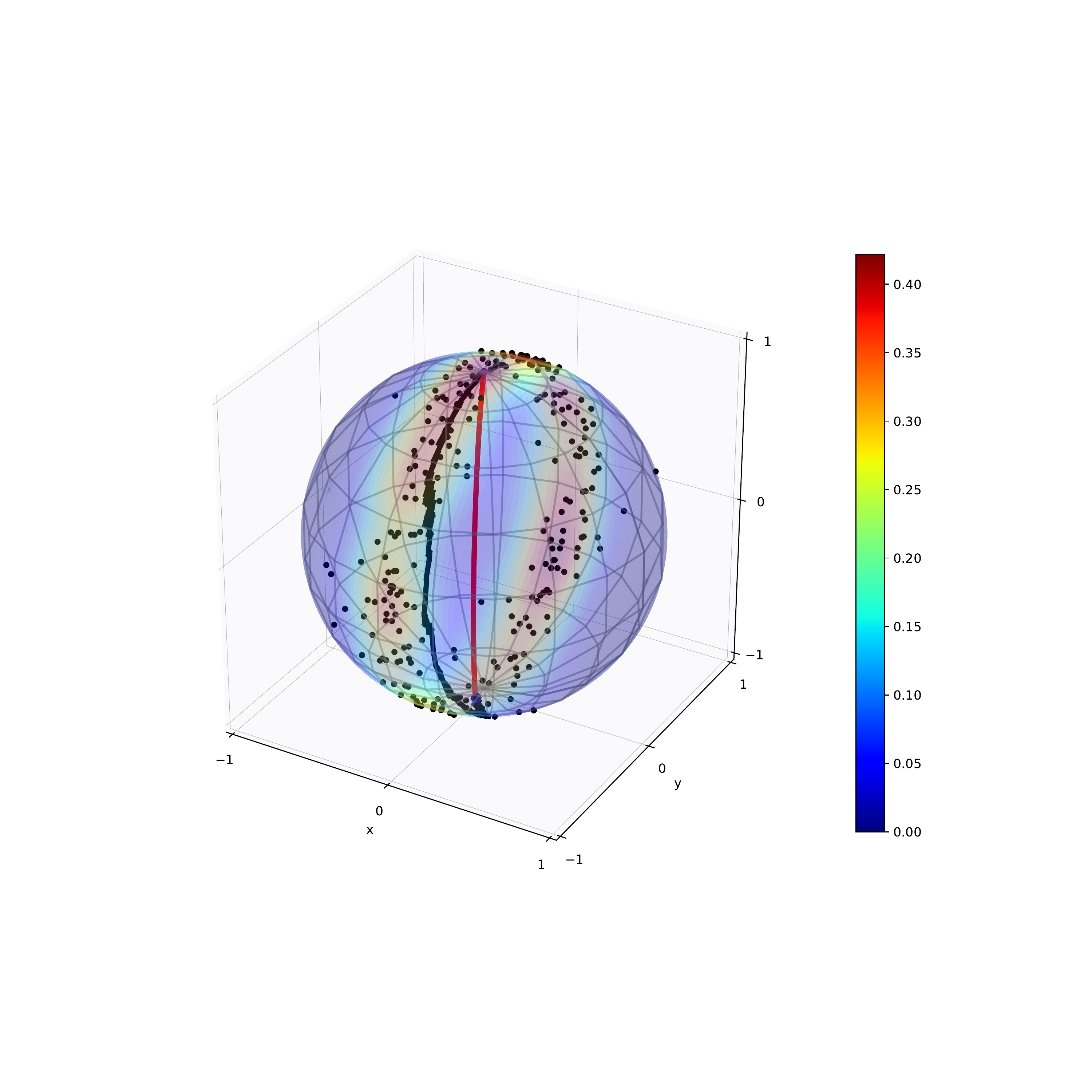}%
    \includegraphics[width=.48\columnwidth, trim = 50 50 50 50,clip]{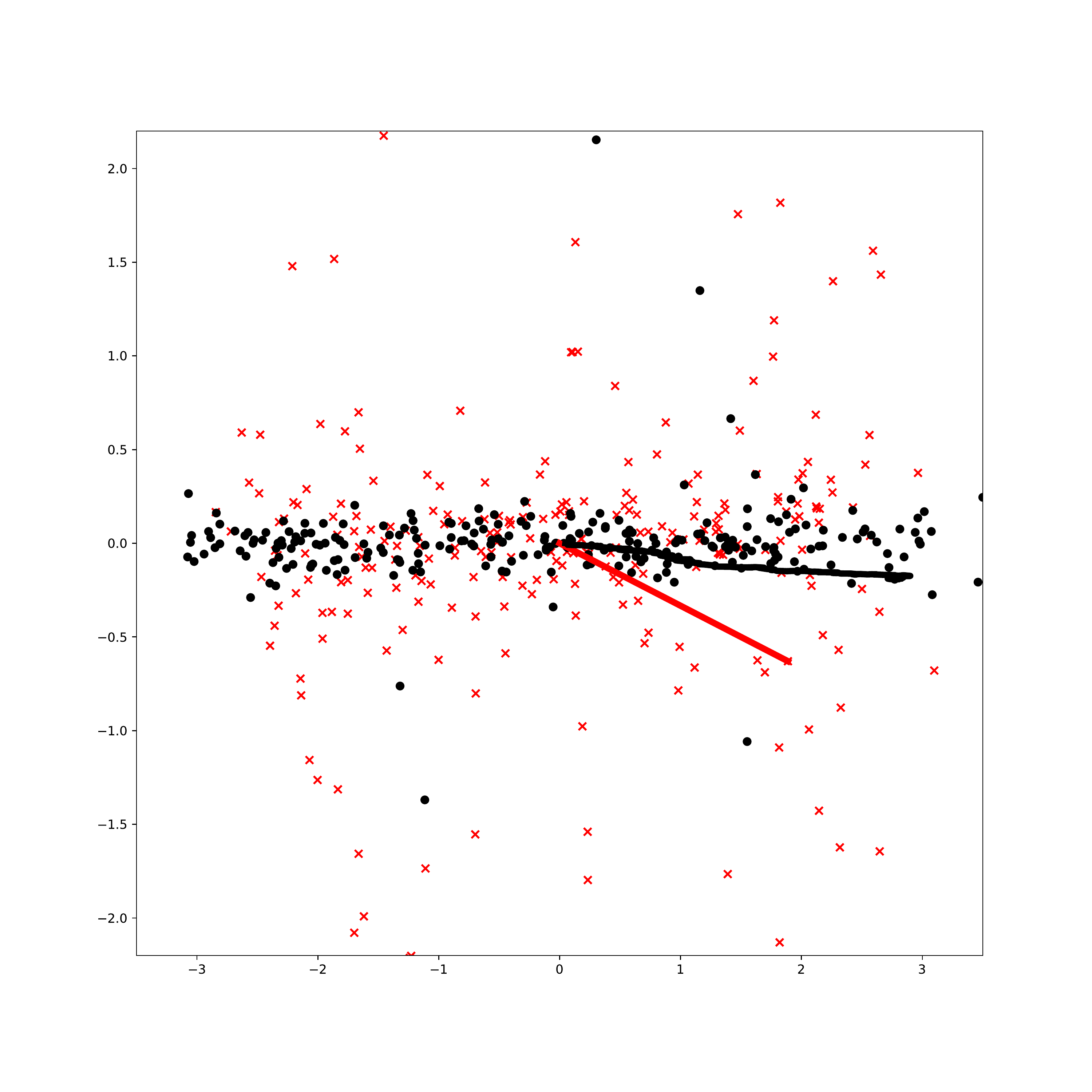}%
    \caption{(left) Samples (black dots) distributed with major mode of variation around a great circle of the sphere and smaller variation orthogonal to the circle. Sphere colored by density of the distribution. (right) Red crosses: The data linearized to the tangent space of the north pole using the Riemannian logarithm map. Because of the curvature of the sphere, variation orthogonal to the great circle is overestimated. This is exemplified by geodesics to data (red straight line/curve in left figure) leaving high-density areas of the data distribution. Black dots: Data linearized to the tangent space using the proposed PCA model. Incorporating the data anisotropy gives a faithful linear view of the data variation. The black curve represents expectation over samples of the latent process conditioned on the same observation as the red curve. The corresponding path is shown on the left figure where it clearly follows the high-density area of the distribution contrary to the geodesic.}
    \label{fig:intro}
\end{figure}

Figure~\ref{fig:intro} visualizes the effect of incorporating anisotropy with the proposed model as compared to generalizing PCA with tangent space linearization. Because of the positive curvature of the sphere, the tangent space linearization overestimates the variance in the second component of the data as geodesic paths may leave high-density areas of the data distribution. In contrast, incorporating anisotropy in the PCA procedure gives a linear view with a faithful representation of the data variation.

The probabilistic construction naturally leads to inference procedures formulated as maximum likelihood or maximum a posteriori fits to data. Using intrinsically defined probability distributions on the manifold thereby avoids some of the complexities that makes non-probabilistic parametric constructions on manifolds inherently complex. For regression, a similar approach has been pursued in \cite{kuhnel_stochastic_2017}. The present paper is partly based on and extends the Oberwolfach abstract \cite{sommer_diffusion_2014}. While the analogy to PPCA is mentioned in \cite{sommer_anisotropic_2015}, the focus of that paper is on defining normal-like distributions and not to generalize PPCA as is the focus here.
\ \\

The paper starts with a short review of PPCA and PPGA before defining the proposed PCA procedure. Constructing the underlying probability model is the subject of the following sections that uses fiber bundle geometry to represent and transport orientation structures over the manifold. We subsequently discuss the proposed PCA procedure in greater depth before outlining inference methods. The paper ends with simple numerical experiments and concluding remarks.

\section{PCA on Manifolds}
Extending Euclidean statistical notions, tools, and inference procedures to the nonlinear manifold situation has been treated in multiple works in recent literature. We focus here on PCA-like statistical analysis of data $y_1,\ldots,y_N$, $y_i\in M$ with $M$ being a nonlinear manifold with a priori known structure, for example arising directly from the data, e.g. angular measurements or position measurements on the surface of the earth, or from modeling constraints. We assume the dimension $d$ of $M$ is finite. Note that the setting is different from manifold learning where the objective is to infer the manifold structure from the data. 

Manifolds lack vector space structure and therefore also a global coordinate system and globally consistent orientations. Instead of inner product structure on Euclidean vectors, the existence of a Riemannian metric $g$ that defines local inner products on infinitesimal variations, vectors in the tangent bundle $TM$, is often assumed. 
While the traditional view of PCA focuses on fitting linear subspaces, we here aim for a probabilistic approach and therefor to generalize the Euclidean probabilistic PCA (PPCA, \cite{tipping_probabilistic_1999}) formulation of PCA to the manifold setting. 
This implies that for the construction in focus in this paper, we mainly need a connection $\nabla$ and a fixed base measure $\mu_0$. If $M$ has a Riemannian metric, $\nabla$ can be the Levi-Civita connection of $g$ and $\mu_0$ the Riemannian volume form $\vol_g$.

We start by outlining Euclidean PPCA and the Riemannian metric based PPGA procedure, before defining the proposed model. We discuss related parametric subspace based methods in section~\ref{sec:parametric}.

\subsection{PCA from a Probabilistic View}
\label{sec:prop_pca}
PPCA interprets PCA as a maximum likelihood fit of the factor model \eqref{eq:linreg} when restricting $W$ to be of rank $k\le d$ and setting the covariance matrix for the noise $\epsilon$ to be diagonal $\sigma^2I$. Because $x$ is assumed normally distributed with unit variance, the marginal distribution of $y$ is normal as well. i.e., PPCA assumes
\begin{equation}
  y|x
  \sim
  N(Wx+m,\sigma^2I)
  \label{eq:y_cond_x}
\end{equation}
with the latent variables $x$ normally distributed
$\mathcal{N}(0,I)$ and i.i.d. isotropic noise $\epsilon\sim\mathcal{N}(0,\sigma^2I)$.
This implies
\begin{equation}
  y
  \sim
  N(m,\Sigma)
  \label{eq:y_marg}
\end{equation}
with $\Sigma=WW^T+\sigma^2I$.

Assuming $\sigma$, $m$ and $W$ are already estimated, the latent variable $x|y_i$ conditioned on the data $y_i$ takes the role of the ordinary principal components of the data $y_i$ in PCA. To get a single data descriptor for $y_i$, one can take the expectation of $x|y_i$ which has the explicit expression $x_i:=E[x|y_i]=(W^TW+\sigma^2I)^{-1}W^T(y_i-m)$. We here loosely denote $x_i$ as principal components for PPCA.

From \eqref{eq:y_marg}, the log-likelihood of the data $y$ is
\begin{equation}
  \mathcal L(y;W,\sigma,m)
  =
  -\frac 12
  (d\ln(2\pi)+\ln|\Sigma|+(y-m)^T\Sigma^{-1}(y-m))
  \label{eq:ppca_likelihood}
\end{equation}
and the maximum likelihood estimate for $W$ is up to rotation given by
$W_{ML}=U_k(\Lambda-\sigma^2I)^{1/2}$,
$\Lambda=\mathrm{diag}(\lambda_1,\ldots,\lambda_k)$ where $U_k$ contains the first $k$ principal eigenvectors of the sample covariance matrix of $y_i$ in the columns, and $\lambda_1,\ldots,\lambda_k$ are corresponding eigenvalues.

For both the ML estimates of $m$ and $W$, and for the principal components $x_i$,
the usual non-probabilistic PCA solution is recovered in the zero noise limit $\sigma^2\rightarrow 0$. 
A similar interpretation of PCA can be found in \cite{roweis_em_1998} where the case $\sigma>0$ is denoted sensible PCA (SPCA).

Turning to the manifold situation, the probabilistic view implies that the fundamental problem in generalizing PPCA is not to define low-dimensional subspaces as sought by the approaches described in section~\ref{sec:parametric} but instead to define a natural generalization of the Euclidean normal distribution to manifolds. PPCA has previously been generalized to manifolds with the probabilistic principal geodesic analysis (PPGA, \cite{zhang_probabilistic_2013}) procedure. The probability model for the data conditioned on the latent variables is for PPGA a Riemannian normal distribution defined via its density 
\begin{equation}
  p(y; m,\tau)
  =
  \frac{1}{C(m,\tau)}
  e^{-\frac\tau 2d_g(m,y)^2}
  \label{eq:Riemannian_pdf}
\end{equation}
with $C(m,\tau)$ a normalization constant and $d_g$ the distance induced by a Riemannian metric $g$ on $M$. This distribution is a function of the squared Riemannian distance to $m$, it is isotropic and closely connected to geodesic distances and least-squares. The latent variables are normally distributed in the linear tangent space $T_mM$ and mapped to the manifold using the Riemannian exponential map $\Exp_m$.

\subsection{An Infinitesimal Probabilistic Model for Manifold PCA}
\label{sec:inf_pca}
We now generalize PPCA using a different probability model. While we follow the PPCA approach of using a maximum likelihood fit of a distribution to data, the distribution here arises from a probability model on infinitesimal steps with covariance on the steps that, when integrated, model the data anisotropy. In contrast to PPGA, the model does not use squared Riemannian distances as in \eqref{eq:Riemannian_pdf} and the Riemannian exponential map.
Instead, we generalize the latent variable model \eqref{eq:linreg} using the anisotropic normal distributions described in section~\ref{sec:anisotropic} to obtain a marginal distribution for the observed data that take the place of the marginal normal distributions in Euclidean PPCA. We here denote this distribution $\mu(m,\Sigma)$ with parameters $m$ and $\Sigma$ for mean and covariance as for the Euclidean normal distribution.

We thus generalize the Euclidean PPCA setting \eqref{eq:y_marg} by assuming the marginal distribution of the data $y$ is 
\begin{equation}
  y\sim \mu(m,\Sigma)
  \label{eq:y_marg_inf}
\end{equation}
with $\Sigma=WW^T+\sigma^I$. Here $W$ has a fixed rank $k$. The distribution $\mu(m,W,\sigma):=\mu(m,\Sigma)$ has a density $p_{\mu(m,W,\sigma)}$ from which we obtain the log-likelihood
\begin{equation}
  \ln\mathcal{L}(y;m,W,\sigma)
  =
  \ln\mathcal{L}(y;\mu(m,W,\sigma))
  =
  \ln p_{\mu(m,W,\sigma)}(y)
  \ .
  \label{eq:likelihood_inf}
\end{equation}
As this likelihood incorporates the curvature of the manifold, it does not have a closed form expression as the Euclidean equivalent \eqref{eq:ppca_likelihood}. However, we will devise a scheme to approximate it by simulation of conditioned bridges of the underlying manifold-valued stochastic process generating $\mu$. 

The generalized PCA model is now up to rotation given by a maximum likelihood estimate $W_{ML}$ using the likelihood \eqref{eq:likelihood_inf}. We write $W=U\Lambda$ with
$\Lambda=\mathrm{diag}(\lambda_1,\ldots,\lambda_k)$ where $U_k$ now contains the first $k$ principal eigenvectors from the model in the columns. $\lambda_1,\ldots,\lambda_k$ are corresponding eigenvalues.

PPCA gets a single low-dimensional data descriptor from the conditional expectation $x_i:=E[x|y_i]$. The generalized model has similar descriptors by conditioning the underlying stochastic process $y_t$ for the response on the data $y_i$ at the observation time $T$. This gives a time-dependent latent variable path describing the data  
\begin{equation}
  \bar{x}_{i,t}=E[x_t|y_T=y_i]
  \ .
\end{equation}
An example of this path is shown in Figure~\ref{fig:intro} where we exemplify the effect of incorporating anisotropy represented by $\lambda_1,\ldots,\lambda_k$ in the model. The time-dependence can furthermore be integrated out to obtain a single descriptor as in PPCA by setting $x_i:=\int_0^Td\bar{x}_{i,t}=\bar{x}_{i,T}$.

The noise $\sigma$ has a similar effect as in PPCA. Due to the infinitesimal nature of the model, the noise influences the underlying stochastic process at each time point $t$.
The anisotropy of the distribution $\mu(m,W,\sigma)$ is represented in the eigenvalues $\lambda_1,\ldots,\lambda_k$ of the matrix $W$. The main complexity is now modelling how this covariance interacts with the curvature of the manifold. Below, we develop the necessary machinery to achieve this and thereby construct the distribution $\mu(m,W,\sigma)$. Further details of the model follows after this in section~\ref{sec:MPPCA}.


\subsection{Parametric Subspace Constructions}
\label{sec:parametric}
We here give a short overview of related non-probabilistic approaches to generalizing PCA to manifolds. Perhaps the most immediate way to handle the lack of coordinate system on manifolds is to use tangent spaces to linearize the manifold and thereby implicitly define a local sense of linear coordinate system on the nonlinear space. This approach is used for generalizations of the principal component analysis procedure in tangent space PCA (tPCA). The principal geodesic analysis (PGA, \cite{fletcher_principal_2004-1}) procedure also uses a tangent space linearization but minimizes the residual distances to the data using manifold distances induced from a Riemannian metric. The central idea in tangent space based procedures is to find a suitable zero-dimensional representation $m$ of the data, often a Frech\'et mean \cite{frechet_les_1948}, and subsequently map the data from the manifold to the linear tangent space $T_m M$. Given a Riemannian structure on $M$, this can be achieved from the geodesic endpoint map $\Exp_m$ and its locally defined inverse $\Log_m$. This construction is however not faithful to the geometry: The effective linearization of the manifold is only locally around $m$ a proper view of the geometry as encoded in the Riemannian metric. The curvature of the manifold will distort the linearized view of the data when significant data mass is observed far from $m$. The fact that the linear view is only one-to-one up to the cut locus of $m$ further emphasizes the approximation in the tangent space linearization.

When using tangent PCA or similar tangent-space based procedures, principal subspaces found as linear subspaces of $T_m M$ are projected to subspaces of $M$ using $\Exp_m$, i.e. as sprays of geodesics originating at $m$. Such subspaces are generally only geodesic at $m$ itself unlike the Euclidean situation where a linear subspace always contains straight lines between all of its points. 
Multiple methods \cite{huckemann_intrinsic_2010,jung_analysis_2012,sommer_horizontal_2013,eltzner_torus_2015,pennec_barycentric_2016} aims at improving this situation by either using particular properties of the data space or by defining other constructions of geometrically natural subspaces. 
Common to these approaches is the explicit construction of low-dimensional subspaces that, focusing on different aspects, are as faithful to the nonlinear geometry as possible.


\section{Fiber Bundle Geometry}
We here review aspects of fiber bundle geometry focusing on the concepts necessary for the intrinsic construction of normal-like probability distributions. We therefore omit many geometric details. More information can for example be found in the papers \cite{sommer_modelling_2017,sommer_evolution_2015}, and the books \cite{hsu_stochastic_2002,kolar_natural_1993}.

The focus is to handle the absence of global orientation fields on the manifold. Indeed, we can choose an ordered basis, a frame, for a single tangent space $T_xM$ providing a reference orientation for vectors in $T_xM$. This however does not provide us with information about vectors in $T_yM$, $x\not=y$. When $M$ is equipped with a connection, the parallel transport along a curve can be used to link the tangent spaces $T_xM$ and $T_yM$. However, the parallel transport is dependent on the curve, and the holonomy of a manifold with non-zero curvature implies that different choices of curves give different parallel transport. This problem is elegantly handled by the Eells-Elworthy-Malliavin construction of Brownian motion that uses a Euclidean martingale and the orthonormal frame bundle to lift the problem to a distribution of orientations over $M$. The outline below is based on and inspired by this idea.

\subsection{The Frame Bundle}
The frame bundle is the set of points $x\in M$ and ordered bases $\nu$ for $T_xM$. For an element $u=(x,\nu)\in FM$, the frame part $\nu$ consists of $d$ basis vectors $\nu_i\in T_xM$. Splitting $u$ in the parts $x$ and $\nu$ technically requires a local trivialization of $FM$. Instead, we let $\pi$ be the projection $FM\rightarrow M$ that just drops the frame from a frame bundle elements and thus sends $u$ to $x$, and we write just $u$ and $u_i$ for the frame and basis vectors. The frame bundle can equivalently be defined as the principal bundle $\GL(\R^d,TM)$ of invertible linear maps between $\R^d$ and the tangent bundle $TM$. An element $u\in \GL(\R^d,TM)$ assigns to a vector $v\in\R^d$ an element $u v\in T_xM$. The $d$ basis vectors $u_i\in T_xM$ in this view appear as the images $ue_1,\ldots,ue_d$ with $e_1,\ldots,e_d$ the standard basis for $\R^d$.

If the manifold is equipped with a connection $\nabla$, each of the basis vectors $u_i$ can be parallel transported along a curve $\gamma$ on $M$ passing $x=\pi(u)$. We write the parallel transport of a vector $v\in T_{\gamma_0}M$ along $\gamma$ as $P_{\gamma,t}(v)$ giving a vector in $T_{\gamma_t}M$. Performing this operation for all $u_i$ gives a transport along $\gamma$ of the entire frame $u$. We can thus lift the parallel transport operation from working on vectors in the tangent bundle $TM$ to transporting frames in $FM$.

The infinitesimal $t\rightarrow0$ limit of the parallel transport of $u$ along $\gamma$ gives an infinitesimal variation in $FM$, i.e. a vector in the tangent bundle $TFM$ of the frame bundle. The span of the tangent vectors arising from such infinitesimal parallel transports, i.e. from choosing curves $\gamma$ on $M$ with different velocities $\dot{\gamma}_0$, defines a linear subbundle of $TFM$ denoted the horizontal subbundle. Another subbundle of $TFM$ is the vertical subbundle $VFM$, and we can write $TFM$ as a direct sum $TFM=HFM\oplus VFM$ thanks to the connection. Elements in the vertical bundle are variations of $u$ that keep $x=\pi(u)$ fixed varying only the frame part in the fiber $\pi^{-1}(x)$ above $x$. Conversely, infinitesimal variations in the horizontal subspace moves $x$ while keeping the frame part of $u$ as fixed as possible as measured by the connection or, equivalently, the parallel transport. $HFM$ variations are thus zero-acceleration as measured by the connection.

An important property of the horizontal bundle $HFM$ is that the pushforward $\pi_*:TFM\to TM$ of the projection $\pi$ is a linear isomorphism when restricted to the horizontal space for a given $u\in FM$, i.e. $\pi_*|_{H_uFM}:H_uFM\to T_{\pi(u)}M$ is invertible. The inverse is called the horizontal lift, here denoted $h_u:T_{\pi(u)}M\to H_uFM$. That is, we can relate vectors in $T_{\pi(u)}M$ and vectors in $H_uFM$ in a one-to-one fashion.
An important consequence, in particular for our purposes, is the fact that the horizontal lift gives a basis of globally defined vector fields $(H_1(u),\ldots,H_d(u))$ for $HFM$. This is very much in contrast to the situation on the base manifold $M$ where topology generally prohibits globally defined non-zero vector fields. We get this basis by, for each basis element $e_i\in\R^d$, using the horizontal lift $u\mapsto h_u(ue_i)$ to get the $HFM$ valued vector field on $FM$ denoted $H_i(u)$. Moreover, if $M$ has a Riemannian metric, the basis is globally orthonormal for $u\in OM$, $OM$ being the subbundle of $FM$ consisting of orthonormal frames, in the sense that $(\pi_*H_1(u),\ldots,\pi_*H_d(u))$ constitutes an orthonormal basis at each point $\pi(u)$. 

\subsection{Sub-Riemannian Structure}
Recall that the density of the Euclidean normal distribution with covariance $\Sigma$ is a function of the weighed quadratic form $x^T\Sigma^{-1}x$. If we let $W$ be a square root $WW^T=\Sigma$, we can write this as $(W^{-1}x)^T(W^{-1}x)$ using the usual $\R^d$ dot product $v^Tv$ of the preimage $W^{-1}x$ regarding $W$ as a linear map $\R^d\to\R^d$. This construction can be naturally extended to give a sub-Riemannian structure on $FM$ that is then, by definition, related to the density of the normal distribution. Because $u\in FM$ can be regarded a linear map $\GL(\R^d,T_{\pi(u)}M)$, we can take $(u^{-1}v)^Tu^{-1}v$ for $v\in T_{\pi(u)}M$. We thus informally regard $u^{-1}$ a square root of the precision matrix $\Sigma^{-1}$ in $T_xM$, or, conversely, $u$ is a square root of the covariance matrix $\Sigma$ that is then a matrix on $\R^d$. To be precise, we define the inner product
\begin{equation}
  \Sigma^{-1}(u)(v,w)=\langle u^{-1}v, u^{-1}w\rangle_{\R^n}
  =(u^{-1}v)^Tu^{-1}v
  \ ,\ v,w\in T_{\pi(u)}M
  \ .
  \label{eq:Sigma}
\end{equation}
The $\Sigma^{-1}$ notation indicates that the inner product should be seen as encoding the precision matrix corresponding to the term $x^T\Sigma^{-1}x$ in the Euclidean normal distribution density. The inner product on $T_{\pi(u)}M$ lifts to an inner product on $T_uFM$
\begin{equation}
  \Sigma^{-1}(u)(v_u,w_u)=\langle u^{-1}(\pi_*(v_u)), u^{-1}(\pi_*(w_u))\rangle_{\R^n}, \qquad v_u,w_u\in T_uFM
  \ .
\end{equation}
Because $\pi_*|_{H_uFM}$is an isomorphism onto $T_{\pi(u)}M$, this product is positive definite on $H_uFM$. However, it degenerates on $VFM$ and therefore does not define a Riemannian structure on $FM$. It does however defines a sub-Riemannian structure. The sub-Riemannian metric can also be viewed as a map $g_\Sigma:TFM^*\rightarrow HFM$ defined by $\xi_u(v_u)=\Sigma^{-1}(u)\left<v_u,g_\Sigma(\xi_u)\right>$ for $\xi_u\in T_uFM^*$. This in addition defines a cometric, also denoted $g_\Sigma$, in the form of an inner product on $T^*FM$ by $g_\Sigma(\xi_u,\eta_u)=\eta_u(g_\Sigma(\xi_u))$, $\xi_u,\eta_u\in T_uFM^*$. Being inverse to the metric which is modeled after precision matrix, the cometric can be seen as encoding covariance.

\subsection{Bundles of Symmetric Positive Definite Tensors}
\label{sec:tensor}
The quadratic form $\Sigma^{-1}$ is an element of the bundle $\Symp M$ of covariant 2-tensors on $M$. The projection $\pi:FM\to M$ can be factored through this bundle giving a map $q$ such that $FM \xrightarrow{\Sigma^{-1}} \Symp M \xrightarrow{q} M$ with $\pi=q\circ\Sigma^{-1}$. It is natural modeling covariance structure using $\Symp M$ since the bundle omits the implicit rotation that a representation of $\Sigma$ by a square root $u\in FM$ imply. Indeed, $\Symp M$ can be viewed as the quotient $FM/O(d)$ where the orthogonal group $O(d)$ acts on the right by $R.u=uR$ for $R\in O(d)$ and $u\in\GL(\R^d,TM)$.

As shown in \cite{sommer_modelling_2017}, the horizontal/vertical splitting of $TFM$ and the sub-Riemannian structure on $FM$ descend to corresponding structures on $\Symp M$. We can therefore work on the two bundles interchangeably in the same way as one shifts between a square root covariance $W$ and the covariance matrix $\Sigma=WW^T$ in Euclidean statistics. As we will see below, the frame bundle supports the development construction for mapping Euclidean semimartingales to the manifold. We therefore often work on $FM$ keeping in mind that the generated covariance structures can be seen as element of the quotient bundle $\Symp$.

\subsection{Development and Stochastic Development}
Let $x_t$ be a Euclidean semimartingale on $\R^d$ defined as the solution to the Stratonovich SDE
\begin{equation}
  dx_t
  =
  b(t,x_t)
  +
  W(t,x_t)
  \circ
  dB_t
  \label{eq:Stratonovich_SDE}
\end{equation}
where $B_t$ is a standard Brownian motion on $\R^k$, $k\le d$, and $\circ$ denotes Stratonovich multiplication. We let $\mathbb P_{x_t}$ denote its law and, when $k=d$, $p_{x_t}(v;x,T)$ denotes the time $T$ transition density of the process evaluated at $v\in\R^d$ when started with initial conditions $x_0=x$. Note that when the drift $b$ is zero and $W$ is a time- and spatially stationary full-rank matrix, $x_T$ is normally distributed with covariance $\Sigma=WW^T$ and density $p_{x_t}(v; x, T)=(2\pi T|\Sigma|)^{-\frac d2}e^{-\frac12(v-x)^T\Sigma^{-1}(v-x)}$. This view of the normal distribution arising as the combined effects of an continuous sequence of infinitesimal random steps with covariance $\Sigma=WW^T$ is particularly well-suited for generalizing to the manifold situation.

We achieve this generalization using the stochastic development construction, see e.g. \cite{hsu_stochastic_2002}. Recall above the existence of a globally defined basis $(H_1,\ldots,H_d)$ for the horizontal bundle $HFM$. This can be used to define an $FM$ valued process from the semimartingale $x_t$ via the SDE
\begin{equation}
  du_t
  =
  H_i(u_t)
  \circ
  dx_t^i
  \ .
  \label{eq:stochastic_development}
\end{equation}
Note the Einstein summation convention implies a summation over the components $dx_t^i$ and the horizontal basis fields. In the deterministic case ($W=0$ in \eqref{eq:Stratonovich_SDE}), the ODE is denoted just development or ``rolling-without-slipping'' due to the fact that the frame represented by a solution $u_t$ is parallel transported, or rolled, along the manifold. This is a consequence of $H_i$ representing infinitesimal parallel transport. In the stochastic case, when $u_0=u$ is an orthonormal frame with respect to a Riemannian metric, i.e. an element of the orthonormal frame bundle $OM$, the construction is the basis for the Eells-Elworthy-Malliavin construction of Brownian motion \cite{elworthy_geometric_1988}. In the following, we denote by $\phi_u(x_t)$ the solution of \eqref{eq:stochastic_development} of a path $x_t\in\R^d$, deterministic or stochastic, started at $u\in FM$. The inverse of $\phi_u$ is denote anti-development.

\subsection{Anisotropic Normal Distributions}
\label{sec:anisotropic}
In \cite{sommer_anisotropic_2015,sommer_modelling_2017}, the stochastic development construction is used with $x_t$ a Euclidean Brownian motion to map from a starting frame $u$ to a distribution $y_T=\pi(u_T)$ on $M$ with density $p_{y_t}(v; u,T)$ with respect to a fixed base measure (e.g. the Riemannian volume form $\vol_g$). The orthonormality condition on $u$ in the Eells-Elworthy-Malliavin construction of Brownian motion is thus relaxed. The result is the anisotropic distribution $y_T$ that has nontrivial covariance in the sense that the infinitesimal stochastic displacements of the process have covariance given by the frame $u$. We denote this distribution $\mu$ below. $\mu$ will take the role as the response distribution when generalizing \eqref{eq:linreg}. The base point $m=\pi(u)$ can be interpreted as the mean of $\mu$, and the frame $u$ itself models the infinitesimal square root covariance. The precision matrix is the inner product $\Sigma^{-1}(u)$ on $T_{\pi(u)}M$. The orientation problem that usually prevents us from defining globally non-zero vector fields on $M$ with special properties, e.g. orthonormality, is thus handled by spreading the orientations stochastically in $FM$ with parallel transport and taking the time $T$ distribution before projecting the resulting distribution to $M$. Figure~\ref{fig:density} shows examples of densities of the generated distributions.
\begin{remark}
  The generated distribution is best viewed as generalizing the linear latent variable model \eqref{eq:linreg} with linear relationship between the covariate and response on an infinitesimal level. Though this linearity is the main focus of the model, the model can be given physical interpretations: For example, with $M=\mathbb S^2$, the horizontal process can describe physical objects on the earth surface that move without drift and with stochastic steps taken with covariance relative to internal gyroscopes. In their movements, keeping zero acceleration of the gyro is exactly parallel translation. At the fixed observation time $T$, $\mu$ describes the distribution of positions of the objects.
\end{remark}

Note that the process $u_t$ is actually a semi-elliptic $FM$-valued Brownian motion with respect to the sub-Riemannian metric $g_\Sigma$ on $FM$. The semi-ellipticity arise because the diffusion is generated only in the subspace $HFM$ of $TFM$. The curvature of $M$ is exactly the non-integrability of the horizontal fields $H_i$, and non-zero curvature therefore implies that the process will diffuse out of the horizontal bundle and generate a larger subspace of $FM$. It does however not satisfy the H\"ormander condition on $TFM$, and the diffusion will not fill all of $FM$.

%

\section{Probabilistic Principal Component Analysis on Manifolds}
\label{sec:MPPCA}

We here provide more detail and precise definitions of the PPCA generalization described in section~\ref{sec:inf_pca}. Consider the map $\phi_{x_t,T}:FM\rightarrow\mathrm{Prob}(M)$ that by stochastic development sends
$u\in FM$ to $\pi(u_T)$ where the $FM$ diffusion $du_t=H_i(u_t)\circ dx_t^i$ is started at time $t=0$ at $u$, and $x_t\in\mathbb{R}^d$ is a Brownian motion. The stopping time $T$ can without loss of generality be assumed $T=1$. Recall from the discussion earlier in the paper that $u$ represents the mean $m=\pi(u)$ and the frame $u$ the square root covariance of the distribution $\mu(u)=\pi(u_T)$. The precision matrix is the inner product $\Sigma^{-1}(u)$ given by $u$. We let $\Gamma\subset\mathrm{Prob}(M)$ be the image of $\phi_{x_t,T}$, i.e. the set of distributions $\mu(u)=\phi_{x_t,T}(u)$ resulting from point-sourced diffusions in $FM$ stopped at time $T$. We then assume the observed data is distributed
according to $\mu(u)\in \Gamma$ so that $y\sim \mu(u)=\pi(u_T)$ for a diffusion
$u_T\in FM$ started at $u$.

Let $\mu_0$ be a fixed measure on $M$, e.g. a Riemannian volume form $\vol_g$. For each distribution $\mu\in\Gamma$, we write $p_\mu$ for the density satisfying $\mu=p_\mu\mu_0\in\Gamma$. We can then define the log-likelihood
\begin{equation}
  \ln\mathcal{L}(y;u)
  =
  \ln\mathcal{L}(y;\mu(u))
  =
  \ln p_{\mu(u)}(y)
  \label{eq:likelihood}
\end{equation}
for a sample $y\in M$. Now for samples $y_1,\ldots,y_N$, let $u_{ML}\in FM$ be a maximum
for $\ln\mathcal{L}(y_1,\ldots,y_N; \mu(u))=\prod_{i=1}^N\ln\mathcal{L}(y_i;u)$. Then $u_{ML}$ contains the parameters of a maximum likelihood fit to the data $y_1,\ldots,y_N$ of the parameters of the model in $u$.

In the PPCA model \eqref{eq:y_marg}, the coefficient matrix $W$ was assumed of rank $k\le d$. A similar rank $k$ model in the nonlinear setting can be constructed by instead of modelling $u$ directly, letting $W$ be an element of the bundle $F^kM$ of rank $k$ linear maps $\R^k\rightarrow TM$. In addition, we need to represent the isotropic iid. noise $\epsilon$ with variance $\sigma^2$. Generators of isotropic noise are elements of the orthonormal frame bundle $OM$ with respect to a Riemannian metric $g$ on $M$, confer the Eells-Elworthy-Malliavin construction of Brownian motion. We denote such an element by $R\in O_{\pi(u)}M$ to emphasize its pure rotation, no scaling nature. We then set
\begin{equation}
  \begin{split}
  &dW_t=H_i(W_t)\circ dx_t^i+\sigma H_i(R_t)\circ d\epsilon_t^i
  \ ,
  \\
  &dR_t=h_{R_t}(\pi_*(dW))
  \end{split}
  \label{eq:uWR}
\end{equation}
and start the processes at $(W,R)$.
Here $\epsilon_t$ is a Brownian motion on $\R^d$ modeling increments of the iid. noise while $x_t$ is now a Brownian motion on $\R^k$ modeling the latent variables. This is a direct extension of the PPCA model \eqref{eq:y_marg} and the latent variable model \eqref{eq:linreg} using the stochastic development construction \eqref{eq:stochastic_development}.  Following the notation of section~\ref{sec:inf_pca} and the full rank case $\mu(u)$ above, we set $\mu(m,W,\sigma):=\pi(W_T)$ where $m=\pi(W)$. The response distribution $\mu(m,W,\sigma)$ is thus the distribution of the base points $\pi(W_T)$, i.e. the distribution of elements of $M$ over which the transported matrix $W_T$ is situated.
The horizontal fields and stochastic development are defined on $F^kM$ in a similar way as on $FM$. However, in the no noise situation $\sigma=0$, the generated distributions $\mu(m,W,\sigma)$ would not have strictly positive density. This is similar to the Euclidean PPCA case in the limit $\sigma\to 0$. This is a consequence of the generated process not being full-rank on $HFM$ if $W\in F^kM$, $k<d$, and the isotropic noise $d\epsilon_t$ is not added.

Note that the system \eqref{eq:uWR} could equivalently be formulated as
\begin{equation}
  \begin{split}
  &dR_t=H_i(R_t)\circ ((Wdx_t)^i+\sigma d\epsilon_t^i)
  \ ,
  \end{split}
  \label{eq:dR}
\end{equation}
i.e. by multiplying the increments $dx_t$ of the Euclidean Brownian motion with a fixed matrix $W$ and adding noise $\epsilon_t$. In practice, with \eqref{eq:uWR}, we need only simulate the $W_t$ evolution in $F^kM$ since the isotropic part $R_t$ can be obtained up to rotation by lifting any $g$ orthonormal basis to an element of $OM$. For high dimensional systems, simulating on $F^k$, $k\ll d$ can be computationally much more tractable than on $FM$ (or $OM$). The system \eqref{eq:dR} lives on $OM$ and therefore does not have a similar reduction property. Finally, the system \eqref{eq:uWR} separates the latent process $x_t\in\R^k$ from the geometric flow of the coefficient matrix $W_t\in F^kM$. This view emphasizes the role of the fiber bundle $F^kM$ in modeling the flow over $M$ of the coefficient matrix while the latent process $x_t$ is Euclidean.

We saw earlier that $u$ defined a sub-Riemannian structure on $FM$. The addition in \eqref{eq:uWR} can also be seen as a sub-Riemannian metric on $FM$ on the form $g_W+\sigma^2\tilde{g}$ \cite{sommer_evolution_2015}. Here $\tilde{g}$ is a lift $\tilde{g}(\xi_u,\eta_u)=g(\pi_*\xi_u,\pi_*\eta_u)$ to $TFM^*$ of the Riemannian metric $g$ on $M$, and $g_W$ is a rank $k$ inner product on $HFM$ defined from the map $W$.

In the following, to simplify notation, we mostly refer to the stochastic process as just $u_t$ without distinguishing between the $FM$ valued full rank version $u_t$ solution to \eqref{eq:stochastic_development} and the low-rank version $W_t$ in $F^kM$ solution to \eqref{eq:uWR}.

\subsection{The Principal Components}
Euclidean Probabilistic PCA reduces the dimensionality of the data by considering the latent variables conditioned on the observed data $x|y_i$. This random variable converges to the principal components as $\sigma^2\rightarrow 0$. In the proposed model \eqref{eq:uWR}, the latent process $x_t$ takes the place of $x$. With non-zero curvature, the latent process cannot directly summarize the observations in single vectors: sample paths $x_t(\omega)$ generating paths $\pi(W_t(\omega))$ hitting the same endpoint $y_i$ on $M$ will in general not have the same endpoint $x_T(\omega)$ in $\R^k$, see Figure~\ref{fig:conditioned_antidevelopment} in the experiments, section~\ref{sec:experiments}. However, we can still consider the conditioned latent variable process $x_t|\pi(W_T)=y_i$. Since
the latent process lives in $\R^k$ where we can take expectation (in contrast to on $M$), we can summarize by the mean of latent sample paths reaching $y_i$:
\begin{equation}
  \bar{x}_{i,t}=E[x_t|\pi(W_T)=y_i]
  \label{eq:principal_components}
\end{equation}
Thus $\bar{x}_{i,t}$ take the role of the latent variables in PPCA.
Note that given the source $W\in F^kM$, the sample paths can be 
equivalently viewed as paths $\pi(W_t(\omega))$ on $M$ or as paths $x_t(\omega)$ in $\mathbb{R}^k$. Examples of mean paths are illustrated in
Figure~\ref{fig:conditioned_samples}. In $\mathbb{R}^k$, the data can be further summarized by 
integrating out the time dependence from $\bar{x}_i$ summarizing the data $y_i$ only in the latent endpoint $x_{i,T}$ with mean $\bar{x}_{i,t}$. The conditioned latent variables in this way provide a Euclideanization of the data similar to those provided by parametric subspace constructions of manifold PCA, section~\ref{sec:parametric}. Because of the process nature of the latent variable, the linearization will be quite different from the linearizations provided by the such methods.

\subsection{Zero Noise Limit}
In the $\sigma\rightarrow 0$ limit, PPCA recovers the original PCA formulation with projections to the latent space that either minimize residual error or maximize variance of the projected data. In the nonlinear case \eqref{eq:uWR}, as $\sigma$ tends to $0$, we get an $FM$ diffusion that progressively concentrates its infinitesimal displacements around $H_{W_t}F^kM$. The short-time asymptotic limit of the likelihood behaves Gaussian-like \cite{sommer_modelling_2017} in the sense 
\begin{equation}
  \lim_{t\to 0}
    2t\log p_{\mu(u),t}(y)=-d_{\Sigma^{-1}}\left(u_0,\pi^{-1}(y)\right)^2
  \label{eq:short_time_asymp}
\end{equation}
with $d_{g_{\Sigma^{-1}}}$ the distance on $FM$ induced by the sub-Riemannian metric $g_\Sigma$.
We can then conjecture that, in the limit, we recover projections to the latent space in a similar sense. If we let $Q(u)$ denote the subspace of $FM$ reachable by horizontal paths starting at $u\in FM$, a natural limit notion of the principal components would be
    \begin{equation}
      \argmin_{\tilde{u}\in Q(u)}
      d_{g_{\Sigma^{-1}}}(\tilde{u},\pi^{-1}(x))^2
      \ .
      \label{eq:limit_principal_component}
    \end{equation}
     We return to the space $Q$ briefly below, and leave the question if the actual $\sigma\to0$ limit of the principal components take a form similar to \eqref{eq:limit_principal_component} to future work.

\subsection{Rotations and Subbundles}
As discussed in section~\ref{sec:tensor}, the above construction is over specified in the sense of an arbitrary rotation being present in the representation $u$ of the covariance similarly to the Euclidean case of specifying covariance with a square root $W$ instead of the actual covariance matrix $\Sigma=WW^T$. We can handle this by quotienting out $O(d)$, instead specifying the construction on $\Symp$, see section~\ref{sec:tensor}. This has however little influence in practice where the rotation implicit in the matrix $u$ can just be ignored.

\subsection{Curvature and Nonintegrability}
While principal subspaces in the Euclidean case are linear subspaces of $\R^d$, the space $Q(W)$ of endpoints of curves starting at $W\in F^kM$, staying horizontal, and generated by the flow equation \eqref{eq:stochastic_development} is not in general a $k$-dimensional submanifold of $F^kM$. The geometric reason is that curvature is equivalent to non-integrability of the horizontal distribution of the vector horizontal fields $H_1,\ldots,H_k$ on $F^kM$, i.e. the $VFM$ valued Lie brackets $[H_i,H_j]$, $1\le i,j\le k$ are non-zero for some $i,j$. Thus, the Frobenius theorem tells us that the span doesn't integrate to a $k$-dimensional submanifold. This can be seen as the key consequence of curvature for PCA like constructions defined via infinitesimal flows. In the present case, we do not need to truncate the non-integrable span to obtain a $k$-dimensional submanifold as is done when e.g. considering geodesic sprays starting at $\pi(u)$. Instead, we simply model the data as being distributed according to the development of horizontal stochastic flows and thus avoid referring to subspaces in the PCA construction. Note that we can still extract principal components as discussed above.

In some cases, we can say more about the structure of $Q(W)$ or $Q(u)$. If $k=d$ and the bracket span $\mathrm{Lie}(H_1,\ldots,H_d\}$ of the horizontal fields has constant dimension for any point at $M$, there exists a subbundle of $FM$ on which $H_1,\ldots,H_d$ satisfies the H\"ormander condition. In this case, the reachable set $Q(u)$ is this subbundle of $FM$, and $\pi(U)$ is a submanifold of $M$. Note that nonzero curvature implies that the dimension $\dim(Q(u))$ is greater than $d$. An example of this case can be seen for the sphere $\mathbb S^2$ where the bracket span has rank $3$ and $Q(u)\simeq O\mathbb S^2$ for any orthonormal $u$. However, if the sphere is deformed to be locally flat in a neighborhood of $\pi(u)$, the rank of the bracket span is lowered to $2$ in this neighborhood and the constant rank condition fails.

\subsection{Extensions}
As noted in \cite{tipping_probabilistic_1999}, the probabilistic formulation has advantages beyond the theoretical insight and the ability to perform estimation with MLE. This includes extension to mixed models where data are assumed distributed according to a sum of multiple latent models of the form \eqref{eq:linreg}, in effect allowing different centers $m_1,\ldots,m_j$ or more complicated shaped distributions. Similar flexibility is present in the manifold situation. The stochastic process $u_t$ or $W_t$ can be started at multiple points $u_1,\ldots,u_j\in FM$ and the resulting densities averaged.

\section{Inference and Predictions}
We here describe two estimation approaches. The first is based on the estimators described in \cite{sommer_anisotropic_2015,sommer_modelling_2017} that use the anisotropically weighted energy of the most probable paths as surrogates for the $\log$ data likelihood. The second approach outlines a Monte Carlo method for estimating transition densities from which the likelihood can be optimized. While the former incorporates anisotropy of the model in the distance $d_{\Sigma^{-1}}$, it uses the short-time asymptotic limit 
\eqref{eq:short_time_asymp} making it only suitable for data with limited variability. The latter method does not employ a similarly approximation. It however includes an expectation of the stochastic process which in practice is approximated by Monte Carlo sampling.

\subsection{Most Probable Paths}
In \cite{sommer_modelling_2017}, the short-time asymptotic limit of $p_{\mu(u),t}(\cdot)$ is used to suggest the estimator
\begin{equation}
\textrm{argmin}_{u \in FM}  
\sum_{i=1}^N 
\left(
d_{\Sigma^{-1}}\left(u,\pi_u^{-1}(y_i)\right)^2
-N\log(det(u)_g)
\right)
\label{eq:mpp_estimator}
\end{equation}
for the maximum likelihood fit of an anisotropic normal distribution to data points $y_1,\ldots,y_N\in M$. The $FM$ distance $d_{\Sigma^{-1}}$ is dependent on $u$, and a minimizer for \eqref{eq:mpp_estimator} can be found by iterative optimization. Because the short-time asymptotic limit is used, the estimator is reasonable for data with limited variation around $\pi(u)$. The $d_{\Sigma^{-1}}$ distances are realized by most probable paths on $FM$ \cite{sommer_anisotropically_2016}, a family of paths that generalizes geodesics when $u$ is not orthonormal.

\subsection{Bridge Simulation}
We generally do not wish to restrict to cases where the data variation is small. As data variation and curvature increases, the estimator \eqref{eq:mpp_estimator} will provide a progressively less precise approximation of the optimal likelihood \eqref{eq:likelihood}. Instead, we here describe a bridge simulation scheme based on the conditioned diffusion bridge simulation method of \cite{delyon_simulation_2006} and the maximum likelihood estimation in \cite{sommer_bridge_2017}.

In \cite{delyon_simulation_2006}, simulation of diffusion processes $x_t\in\R^d$ given by the It\^o SDE
\begin{equation}
  dx_t
  =
  b(t,x_t)dt
  +
  W(t,x_t)dB_t
\end{equation}
conditioned on hitting a point $v\in\R^d$ at time $T$
is considered based on the idea of adding a drift term that guides the diffusion towards the target $v$. The resulting modified SDE takes the form
\begin{equation}
  d\tilde{x}_t
  =
  b(t,\tilde{x}_t)dt
  -
  \frac{\tilde{x}_t-v}{T-t}dt
  +
  W(t,\tilde{x}_t)dB_t
  \ .
\end{equation}
Under reasonable assumptions, including that $W$ is invertible for all $t,x$, \cite{delyon_simulation_2006} shows that $E_{x_t|v}[f(x_t)]=E_{\tilde{x}_t}[f(\tilde{x}_t)\varphi(\tilde{x}_t)]$ for measurable maps $f$ on $W([0,T],\R^d)$. Here $\varphi$ is a correction factor that takes into account the difference of the laws of the process $x_t|v$ that is conditioned on hitting $v$ at time $T$, and the modified process $\tilde{x}_t$. Note that $\tilde{x}_t$ by construction will hit $v$ a.s. The density of $x_T$ can be recovered from this construction as
\begin{equation}
    p(v; x,T)
    =
    \left(\frac{\left|\Sigma^{-1}(v)\right|}{2\pi T}\right)^{\frac d2}
    e^{-\frac{\|W(x)^{-1}(x-v)\|^2}{2T}}
    \mathbb E_{\tilde{x}_t}
    [\varphi(\tilde{x}_t)]
    \label{eq:density_Euclidean}
\end{equation}
with $\Sigma=WW^T$.

We now suggest to use a similar approach for estimating the likelihood of the data under the proposed manifold PPCA model. Because we do not have a diffusion process with invertible diffusion field $W$ as above (the process $u_t$ (or $W_t$) is only semi-elliptic on $FM$), we will not here give a rigorous argument for the convergence of the procedure. We will instead sketch an approach that uses the fact that the data in the model is only observed at $M$ while the process $u_t$ lives in $FM$. This situation is related to the case of partial observations treated in \cite{marchand_conditioning_2011}, see also the semi-elliptic phase-space flows in \cite{arnaudon_geometric_2018}. See also \cite{sommer_diffusion_2018} for details on the construction of guided bridge simulation schemes on nonlinear manifolds.

We assume existence of a chart that covers $M$ except for a set of measure zero, and we use this to write the process in coordinates. The equations below are coordinate expressions in this chart. In particular, the difference $x_t-v$ to the target point is a coordinate difference.

The idea is now to extend the coordinates on $M$ to coordinates on $FM$ as in \cite{mok_differential_1978} and write the $u_t$ diffusion as a process in coordinates. The coordinates imply a trivialization of $FM$ so we can write $u_t=(x_t,\nu_t)$ with $x_t\in M$ and $\nu_t$ frames. We write the system in short form as
\begin{equation}
  \begin{pmatrix}
    dx_t\\
    d\nu_t
  \end{pmatrix}
  =
  \begin{pmatrix}
    b_x\\
    b_\nu
  \end{pmatrix}
  dt
  +
  \begin{pmatrix}
    W_x\\
    W_\nu
  \end{pmatrix}
  \circ
  dW
  \ .
\end{equation}
The $W_x$ part is actually just $\nu$ by construction of the process. We then make a modified process $\tilde{u}_t=(\tilde{x}_t,\tilde{\nu}_t)$
\begin{equation}
  \begin{pmatrix}
    d\tilde{x}_t\\
    d\tilde{\nu}_t
  \end{pmatrix}
  =
  \begin{pmatrix}
    b_x\\
    b_\nu
  \end{pmatrix}
  dt
  -
  \begin{pmatrix}
    \tilde{\nu}\\
    W_\nu
  \end{pmatrix}
  \frac{\tilde{\nu}^{-1}(\tilde{x}_t-v)}{T-t}\\
  dt
  +
  \begin{pmatrix}
    \tilde{\nu}\\
    W_\nu
  \end{pmatrix}
  \circ
  dW
  \ .
  \label{eq:tildext}
\end{equation}
Intuitively, the use of $\tilde{\nu}^{-1}$ in the drift term produces a correction that after multiplication on $(\tilde{\nu},W_\nu)^T$ points in the direction $\tilde{x}_t-v$ on $M$ while staying horizontal on $TFM$.

Without arguing for convergence here, we aim for $\tilde{x}_t$ to hit $v$ at time $T$ because of the added drift term. We then find the correction term $\varphi$ as in \cite{delyon_simulation_2006}, and arrive at the expression
\begin{equation}
    p(v; u,T)
    =
    \left(2\pi T\left|u\right|_g^2\right)^{-\frac d2}
    e^{-\frac{\|u^{-1}(x-v)\|^2}{2T}}
    \mathbb E_{\tilde{u}_t}
    [\varphi(\tilde{u}_t)]
    \label{eq:density_lim_E}
\end{equation}
for the density of the process with correction term $\varphi$.

We can now write the density with respect to $\mu_0$, e.g. $\vol_g$ with $g$ a Riemannian metric, and sample from $E_{\tilde{u}_t}[\varphi(\tilde{u}_t)]$ with a Monte Carlo scheme. We do this with Hamiltonian updates to keep the acceptance rate high. We can then optimize for $u=(x,\nu)$ either by directly taking gradients with respect to $u$ of the sample approximation of the likelihood, or by an EM-approach.

\section{Experiments}
\label{sec:experiments}
We aim here to visualize the effect of the method and the influence of curvature on two low-dimensional manifolds, the sphere $\mathbb S^2$ and a non-spherical ellipsoid. While curvature effects are visible in both cases, the non-symmetrical nature of the ellipsoid emphasizes the differences to the Euclidean situation. For both manifolds, we illustrate samples from the model with fixed mean and covariance encoded in the frame bundle element $u$. We then for optimal $u$ illustrate how the non-linearity affects the principal components \eqref{eq:principal_components}. After this, we illustrate iterations of a direct optimization of the approximate data likelihood from the density expression \eqref{eq:density_lim_E} using Monte Carlo sampling.
\begin{figure}[ht]
    \centering
    \includegraphics[width=.45\columnwidth, trim = 130 150 130 150,clip]{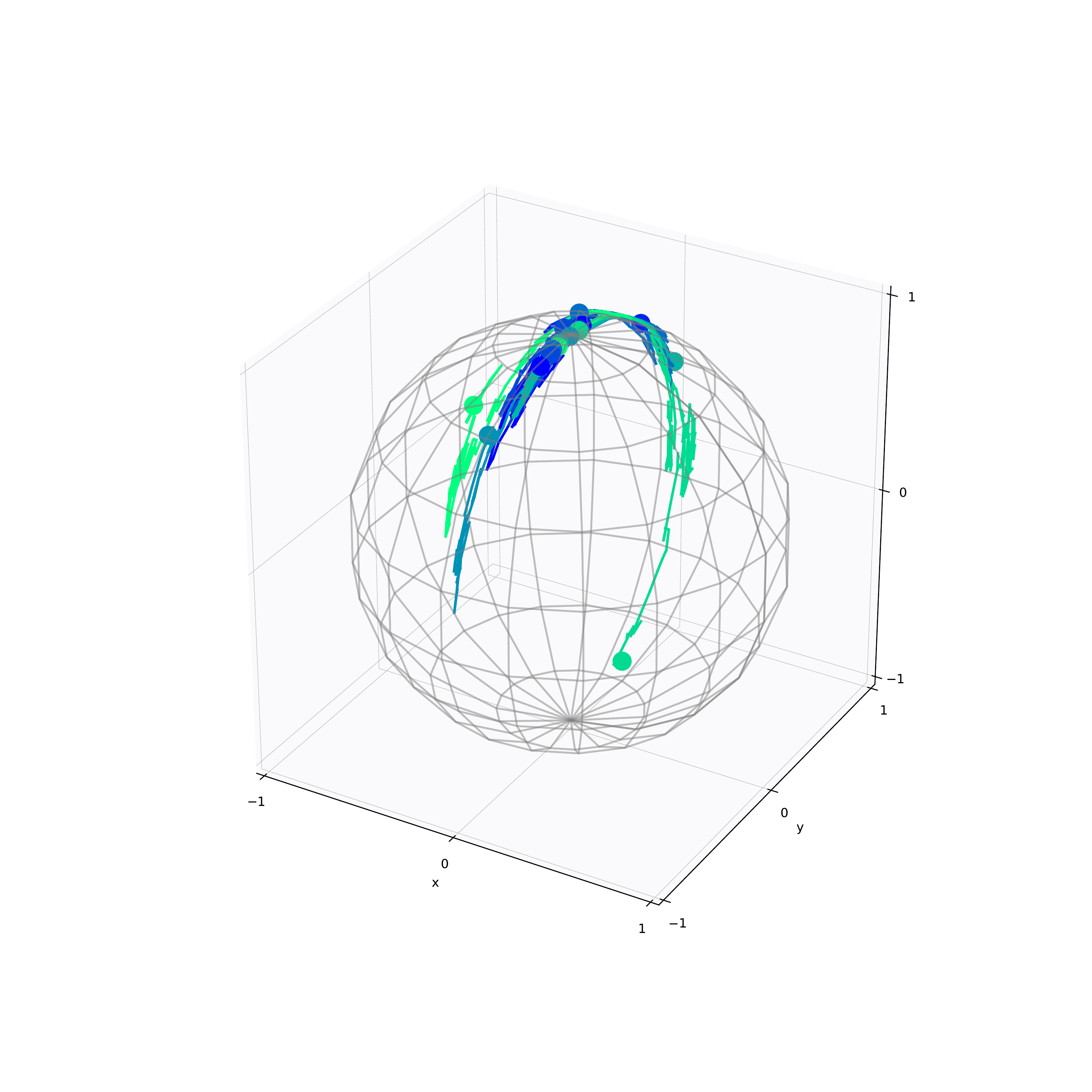}%
    \hspace{.2cm}
    \includegraphics[width=.47\columnwidth, trim = 130 150 90 150,clip]{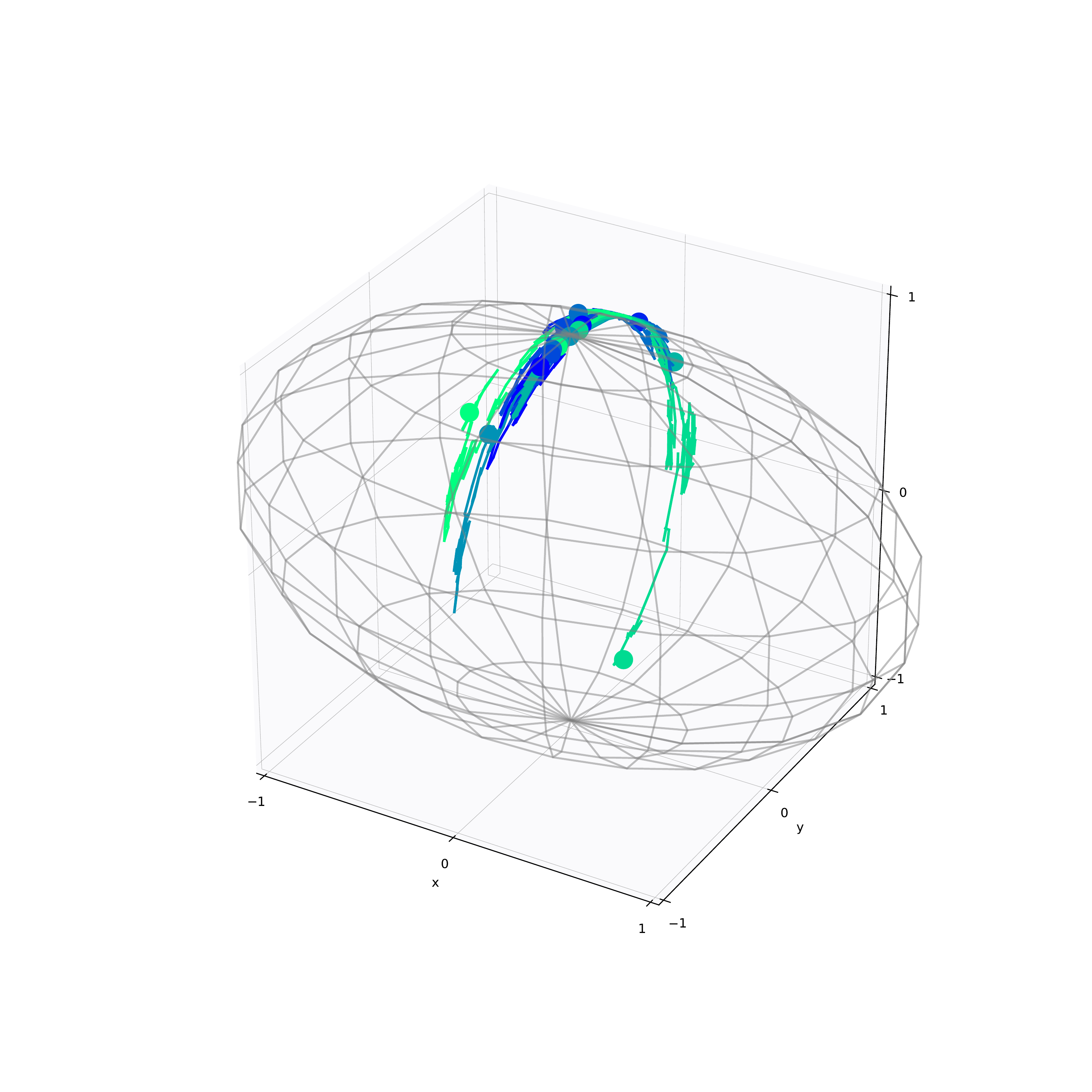}%
    \caption{Samples with corresponding trajectories on the sphere $\mathbb{S}^2$ and an ellipsoid. The variance is 1 in the axis of major variation, and noise with variance $\sigma=.1$ is added in the orthogonal direction.}
    \label{fig:samples_trajectory}
\end{figure}
\begin{figure}[ht]
    \centering
    \includegraphics[width=.48\columnwidth, trim = 130 150 130 150,clip]{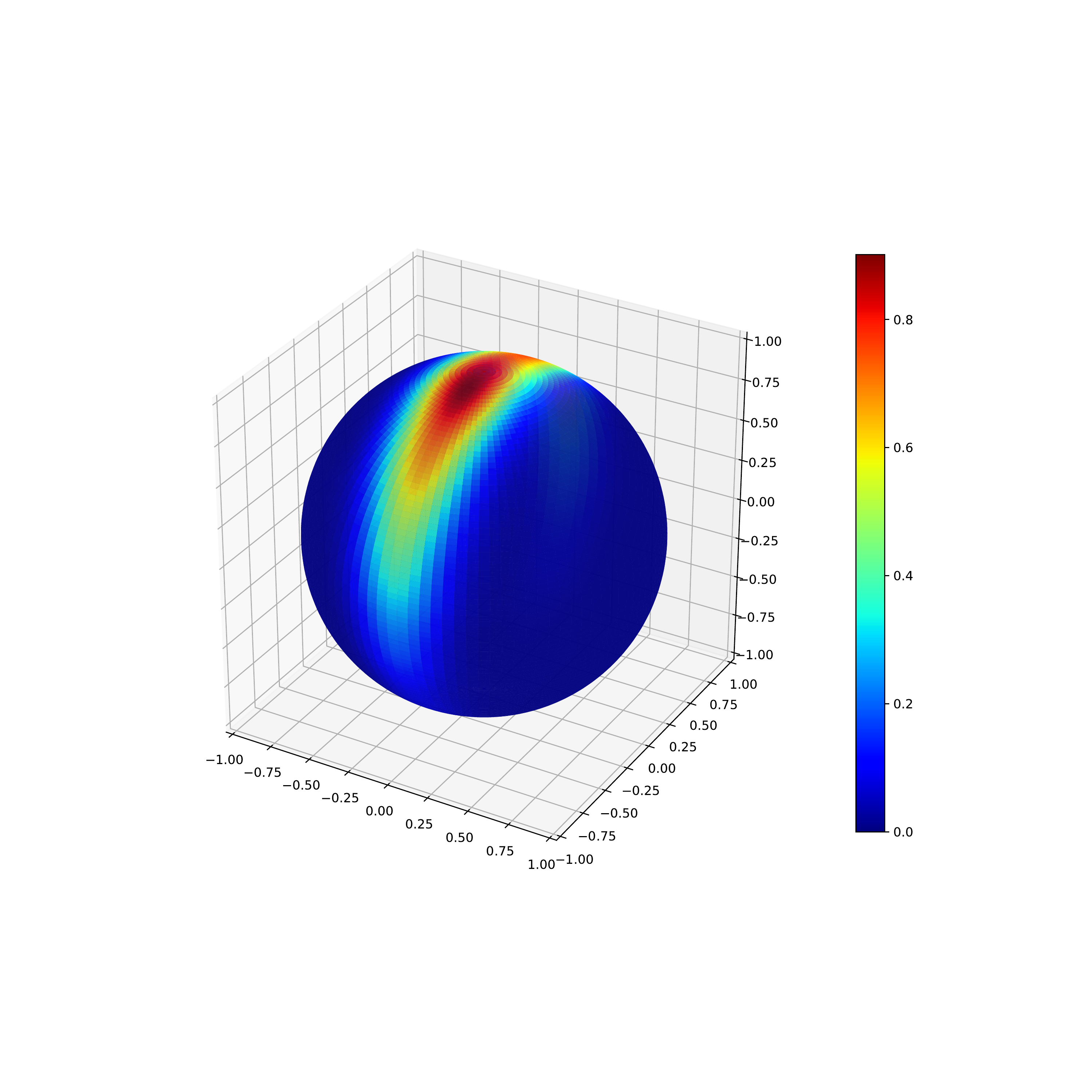}%
    \hspace{.2cm}
    \includegraphics[width=.48\columnwidth, trim = 130 150 130 150,clip]{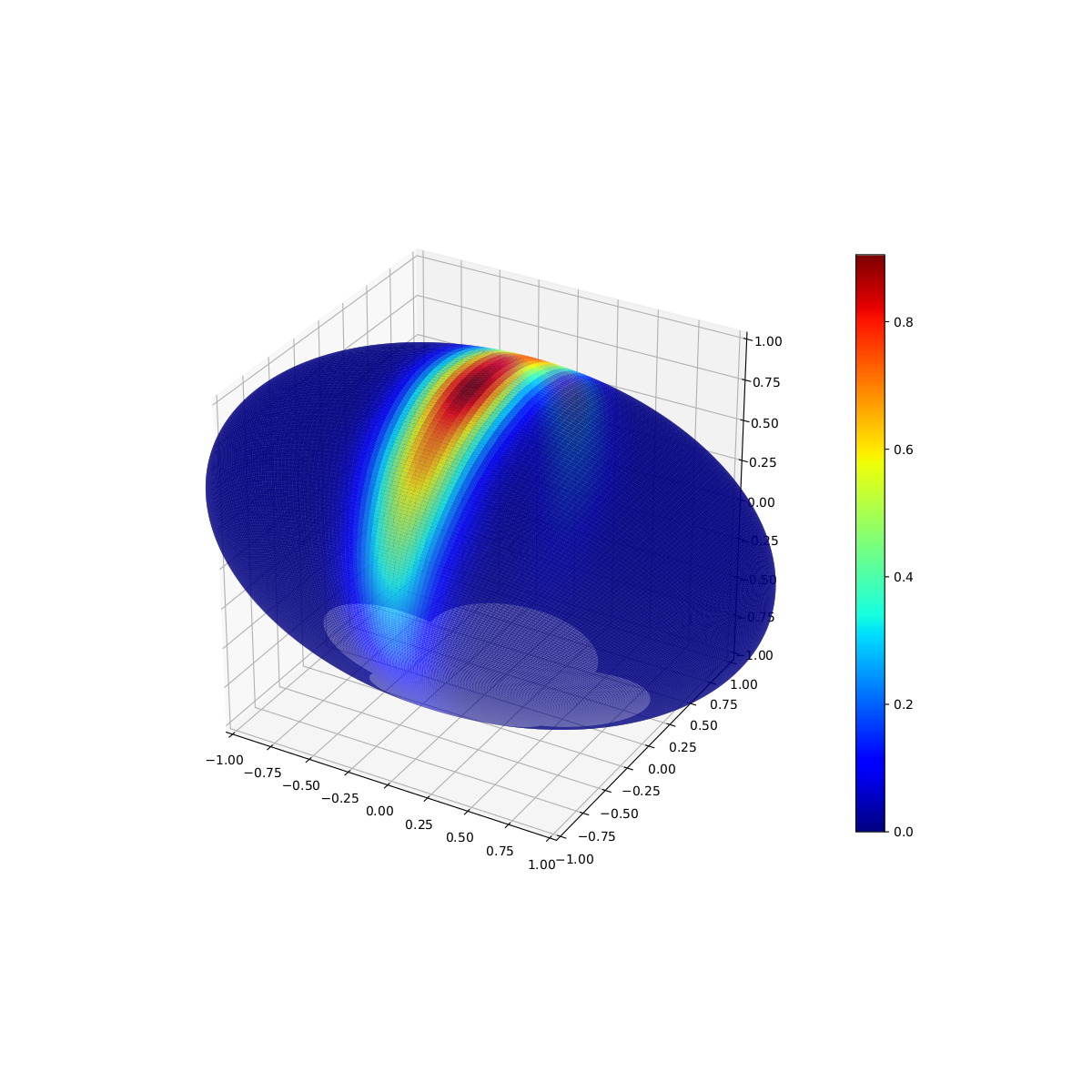}%
    \caption{Density plots on the generated distribution on both surfaces.}
    \label{fig:density}
\end{figure}

The experiments are performed using the differential geometry library Theano Geometry\footnote{\url{https://bitbucket.com/stefansommer/theanogeometry}} that is based on the Theano  framework \cite{the_theano_development_team_theano:_2016} for symbolic expression, automatic differentiation, and subsequent numerical evaluation. See also \cite{kuhnel_differential_2017} for an extended description of the use of automatic differentiation for differential geometric and nonlinear statistical computations. Sampling from the likelihood expression \eqref{eq:density_lim_E} with Hamiltonian updates coupled with gradients for $u$ involves very complex expressions with high order derivatives that would be practically infeasible to derive by hand. Fortunately, the use of automatic differentiation removes this complexity.

\subsection{Density and Forward Sampling}
Figure~\ref{fig:samples_trajectory} shows samples from the probability model on the sphere $\mathbb S^2$ and the ellipsoid with variance $1$ in one axis, and noise with variance $\sigma=.1$ in the orthogonal axis corresponding to the model \eqref{eq:uWR}. The starting point of the diffusion $\pi(u)$ corresponding to the mean is on both surfaces the north pole. The trajectories of the anisotropic process leading to the generated samples are visualized along with the endpoints. Figure~\ref{fig:density} shows the corresponding density on both surfaces.
\begin{figure}[t]
    \centering
    \includegraphics[width=.45\columnwidth, trim = 130 150 90 150,clip]{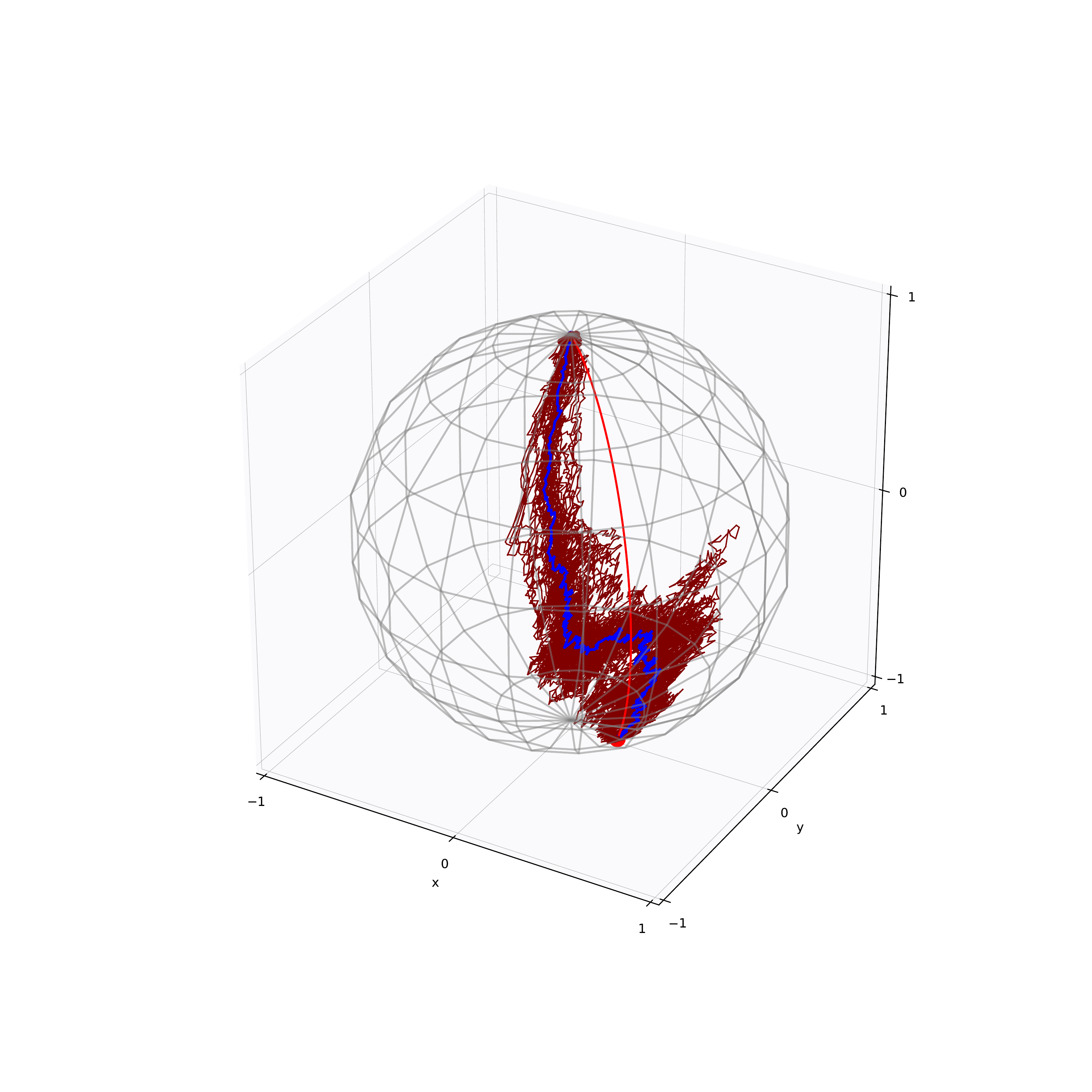}%
    \hspace{.2cm}
    \includegraphics[width=.47\columnwidth, trim = 130 150 90 150,clip]{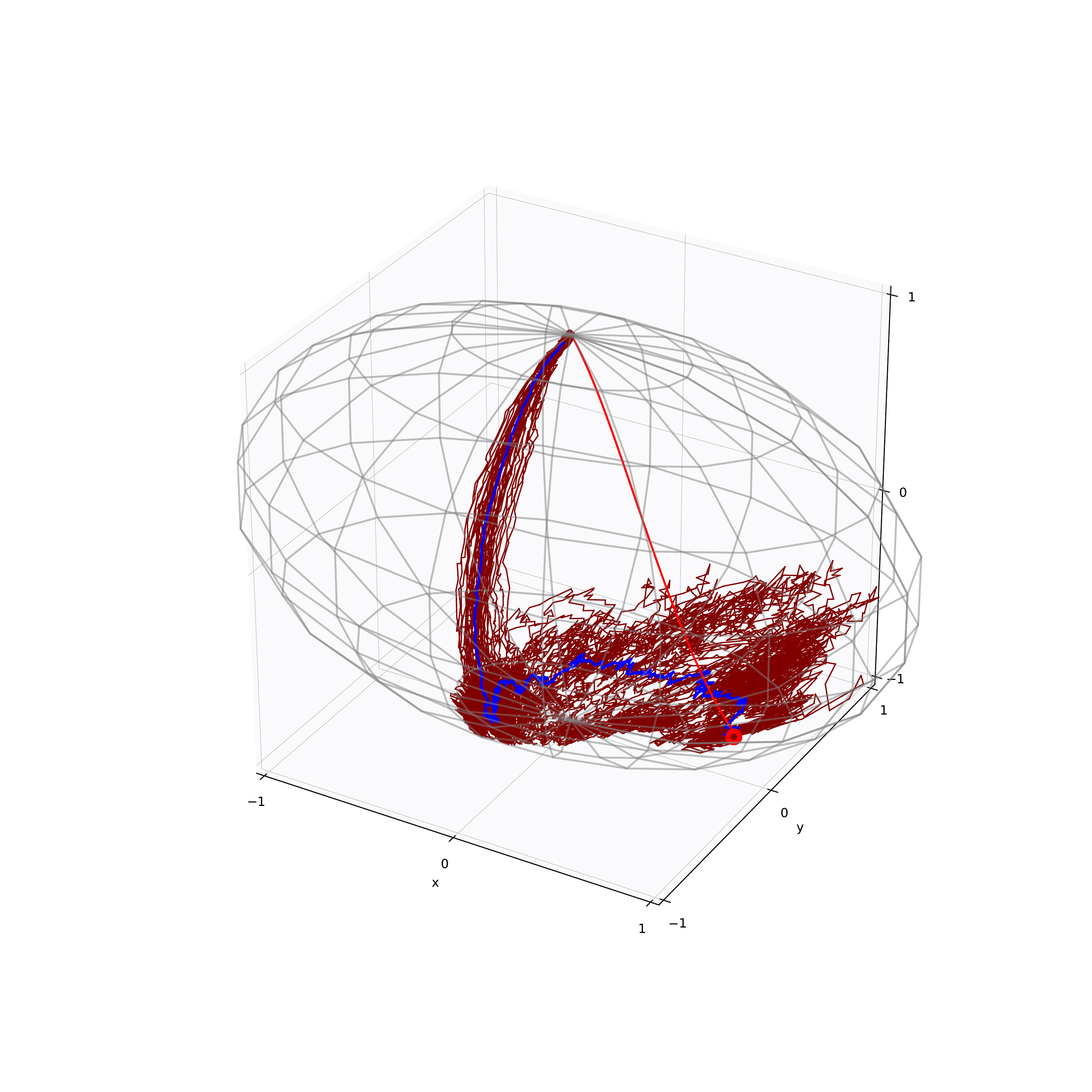}%
    \caption{Samples from the process \eqref{eq:tildext} conditioned on hitting the point $v$. The mean path (blue) plotted deviates from a geodesic (red) to $v$ because of the coupling between the curvature and the anisotropic covariance.}
    \label{fig:conditioned_samples}
\end{figure}
\begin{figure}[t]
    \centering
    \includegraphics[width=.47\columnwidth, trim = 50 50 50 50,clip]{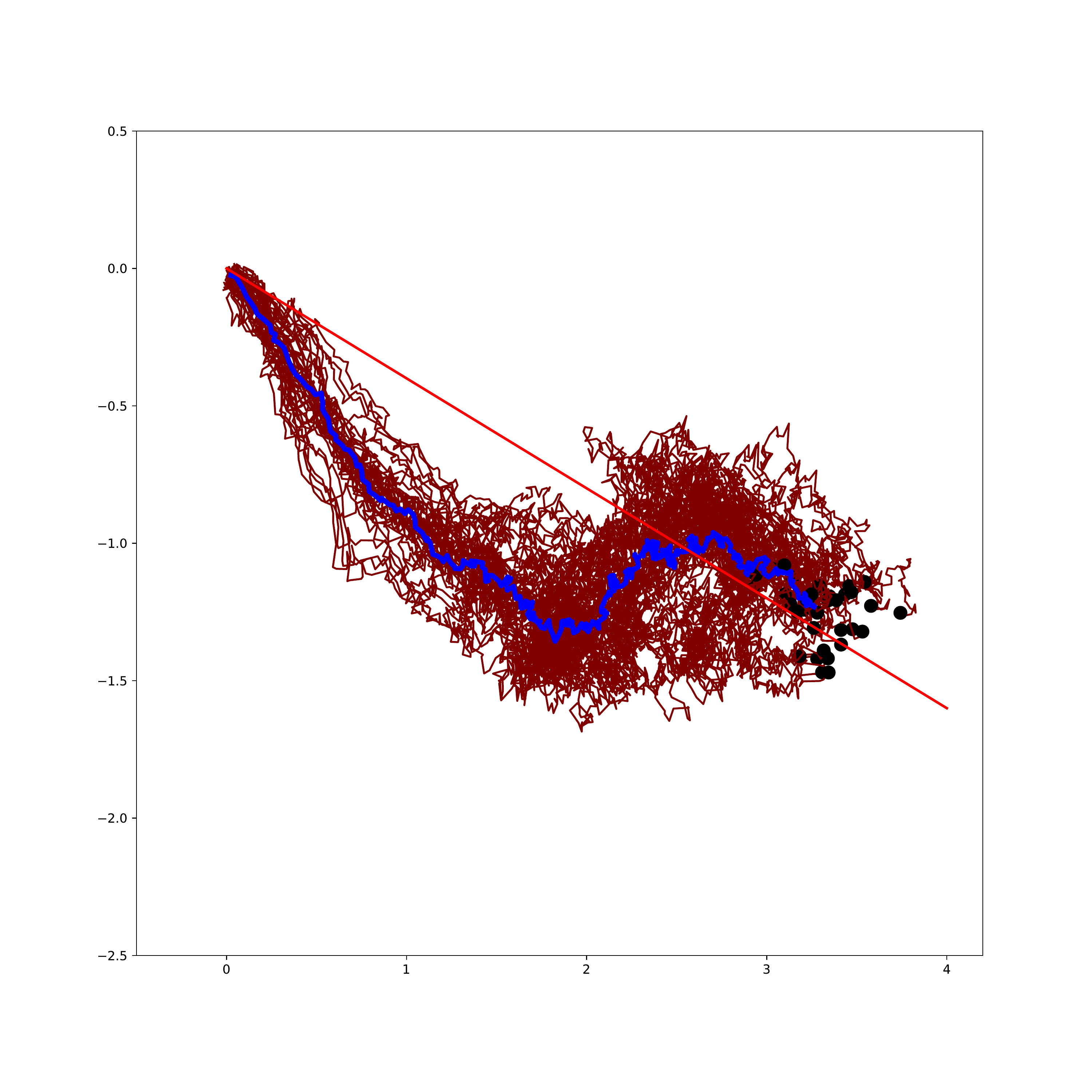}%
    \hspace{.2cm}
    \includegraphics[width=.47\columnwidth, trim = 50 50 50 50,clip]{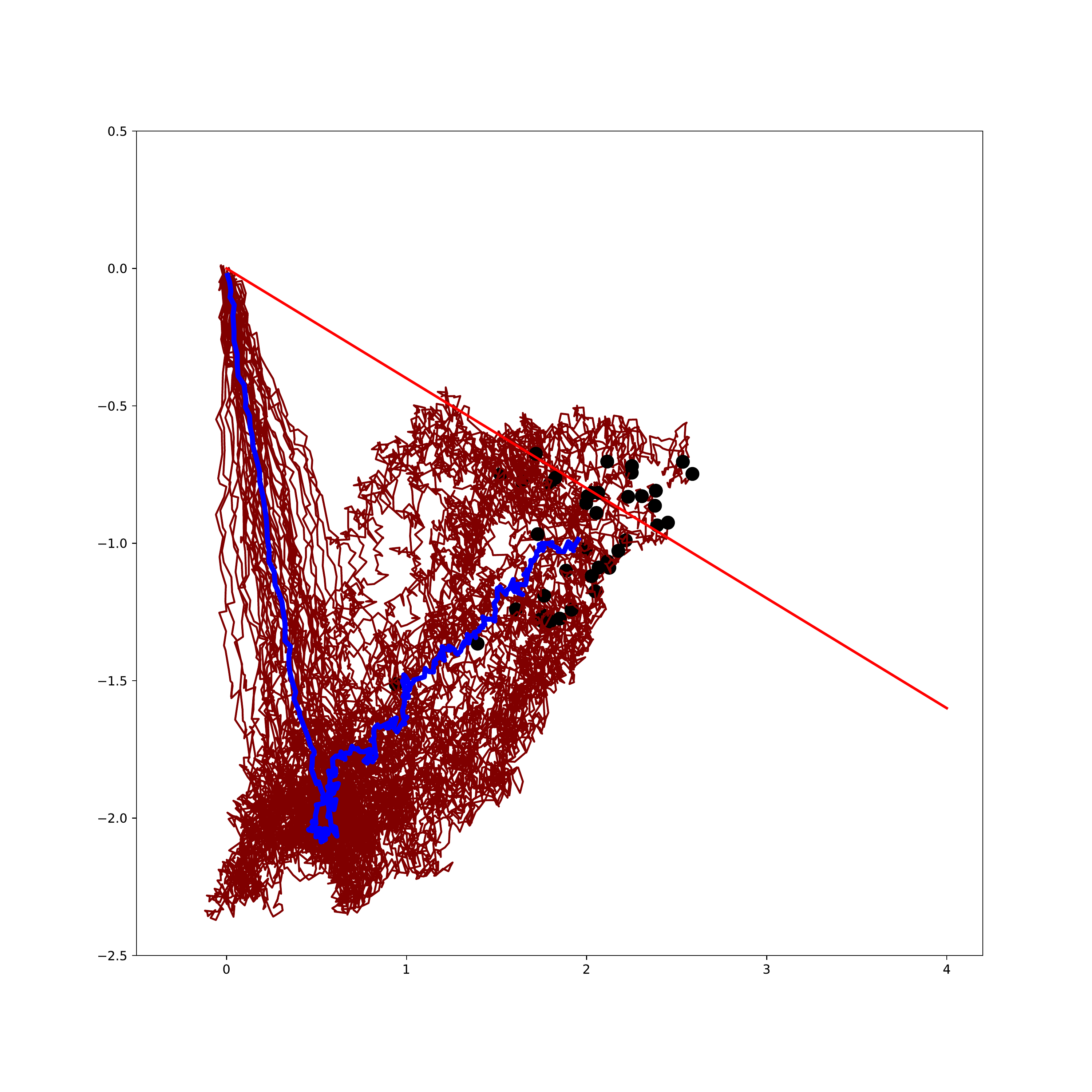}%
    \caption{Samples from the latent process $x_t$ corresponding to the samples in Figure~\ref{fig:conditioned_samples}. The mean latent path is plotted in blue. Even though the process is conditioned on hitting $v$, the endpoints (black) of the latent path deviates. The mean path is not straight and therefore does not correspond to a geodesic on the surfaces.}
    \label{fig:conditioned_antidevelopment}
\end{figure}
\begin{figure}[t]
    \centering
    \includegraphics[width=.47\columnwidth, trim = 50 50 50 50,clip]{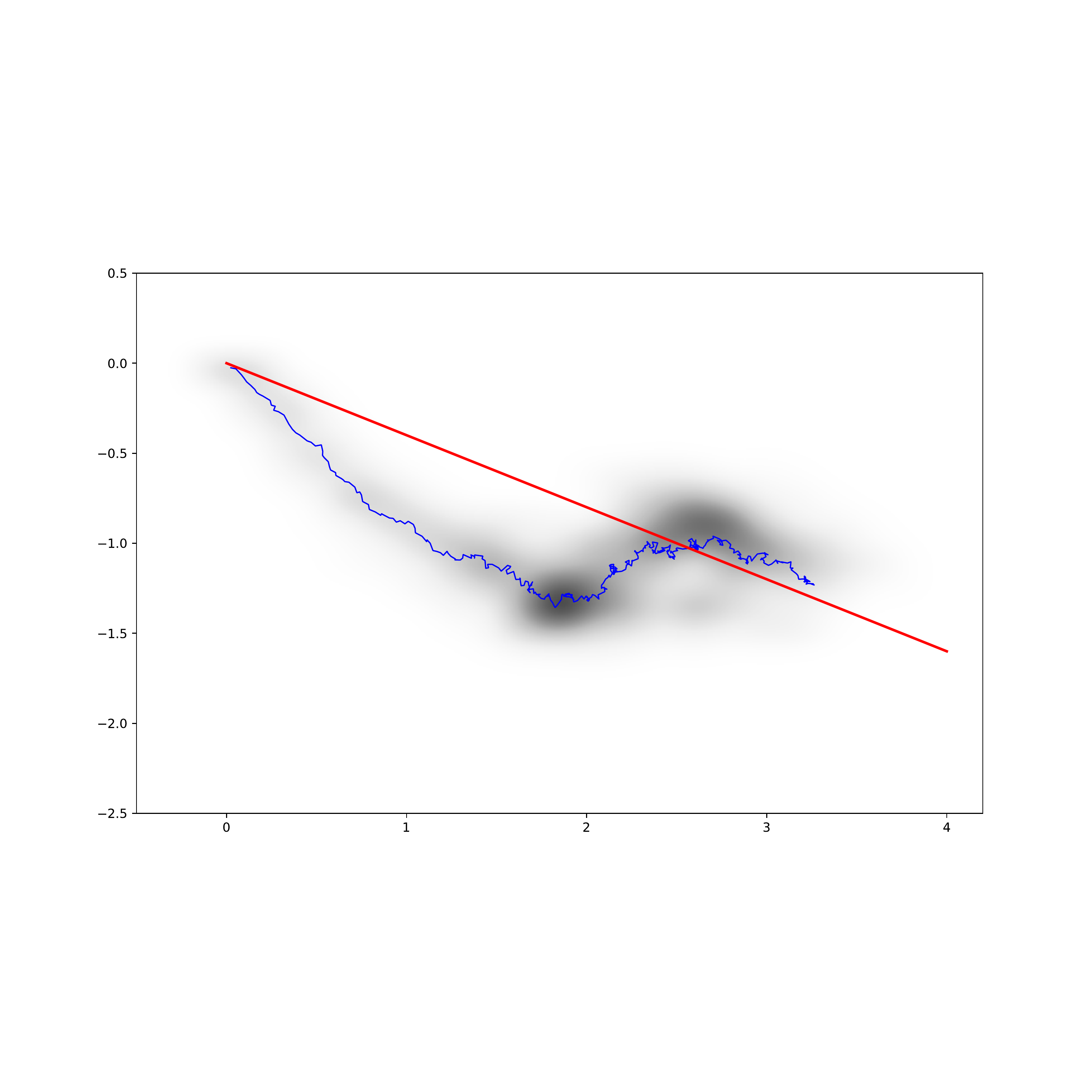}%
    \hspace{.2cm}
    \includegraphics[width=.47\columnwidth, trim = 50 50 50 50,clip]{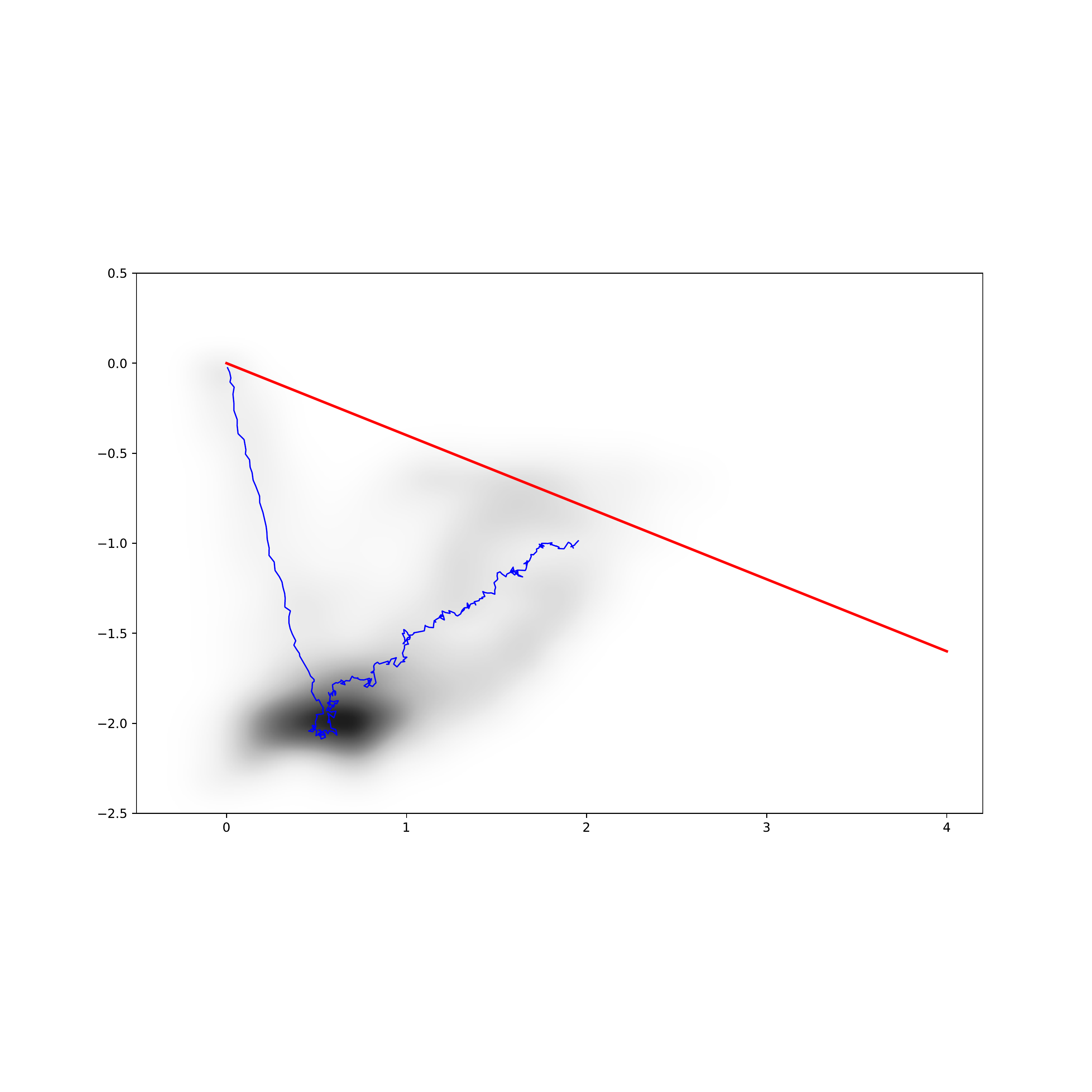}%
    \caption{Density plot of the trajectories in Figure~\ref{fig:conditioned_antidevelopment}. A straight line (red) corresponding to a geodesic on the surface from the north pole to $v$ is plotted for comparison with the mean path (blue solid).}
    \label{fig:conditioned_antidevelopment_density}
\end{figure}
\begin{figure}[h]
    \centering
    \includegraphics[width=.47\columnwidth, trim = 130 150 90 150,clip]{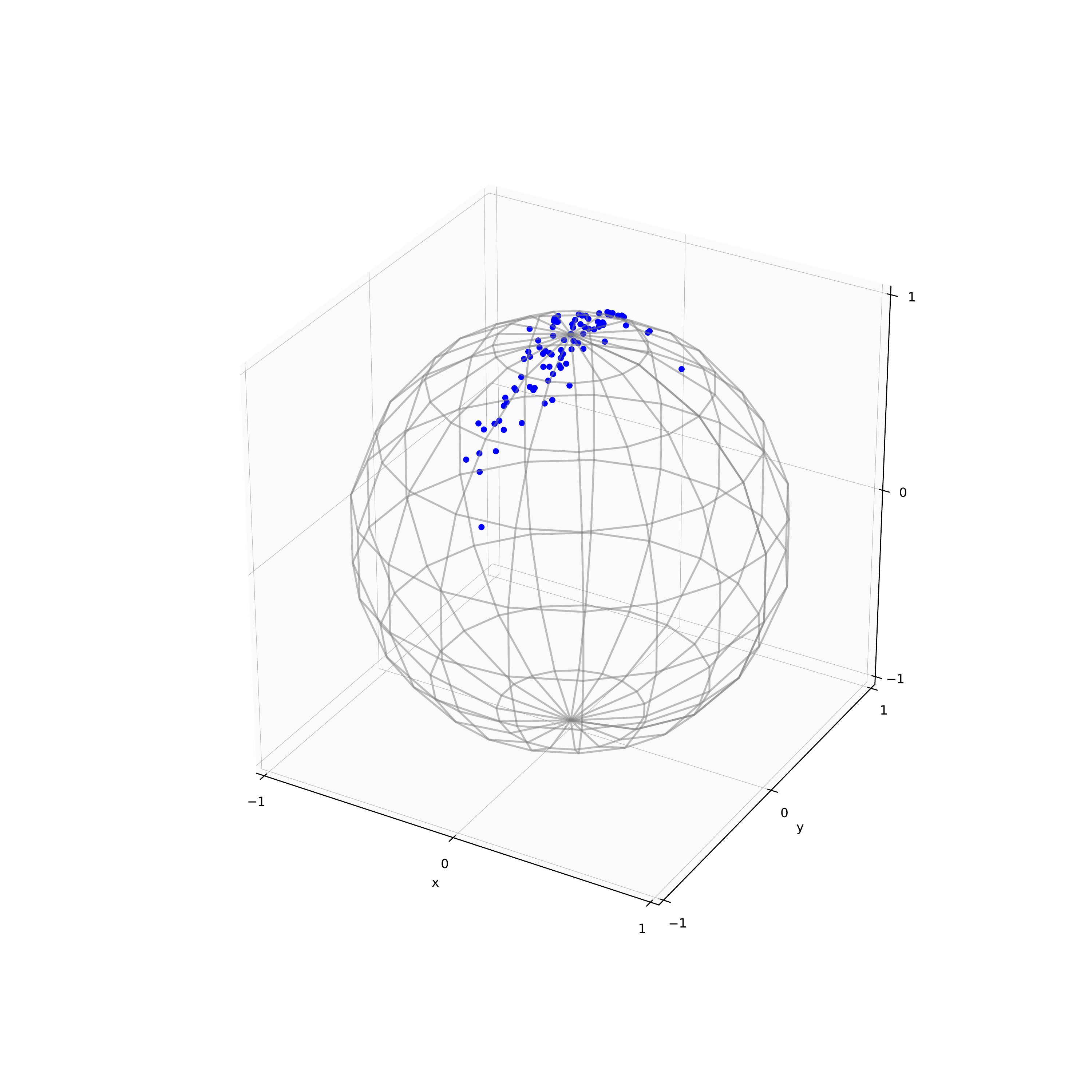}%
    \caption{Samples used for the ML estimation in Figure~\ref{fig:likelihood}.}
    \label{fig:samples_mle}
\end{figure}
\begin{figure}[ht]
    \centering
    \includegraphics[width=.45\columnwidth, trim = 50 50 50 50,clip]{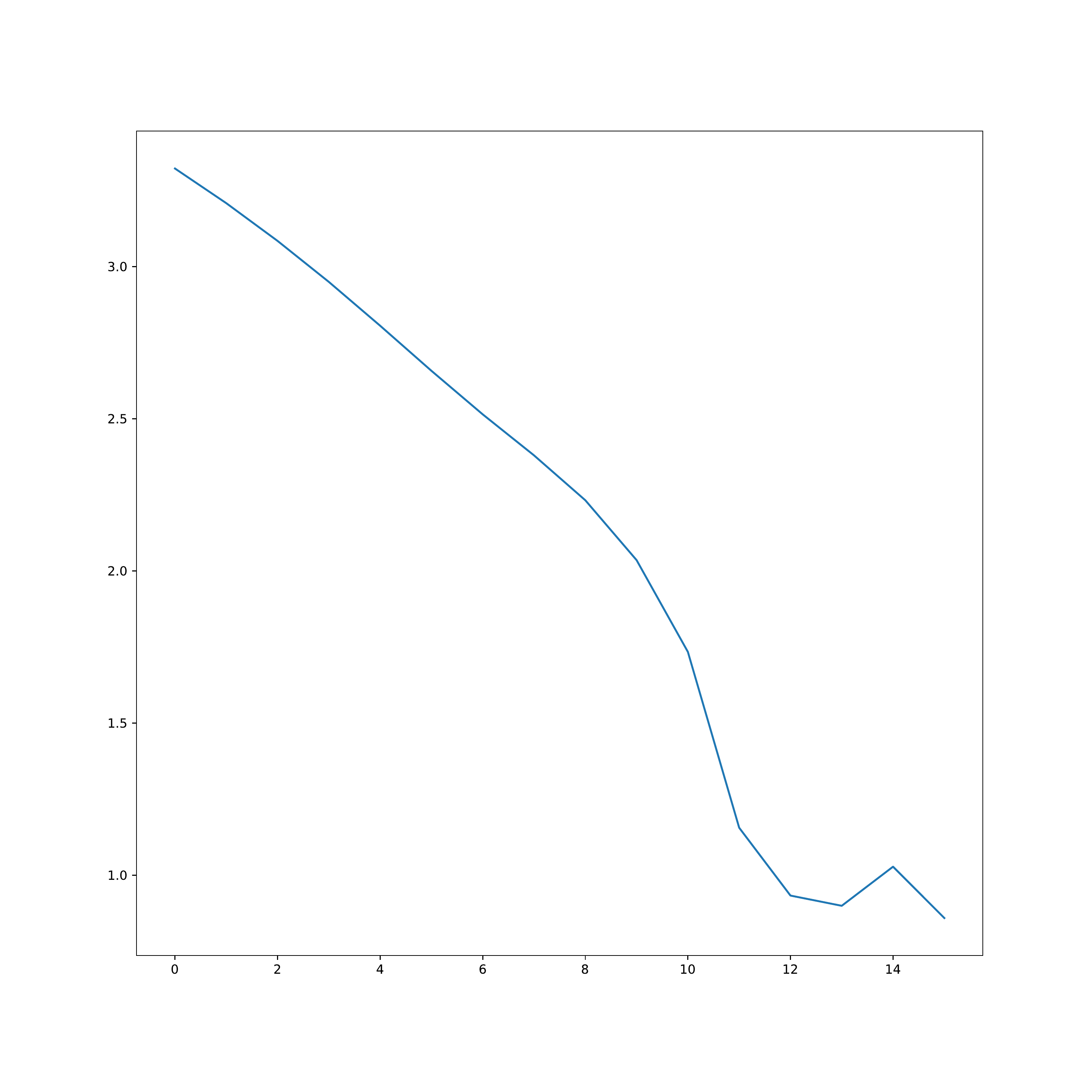}%
    \hspace{.2cm}
    \includegraphics[width=.45\columnwidth, trim = 50 50 50 50,clip]{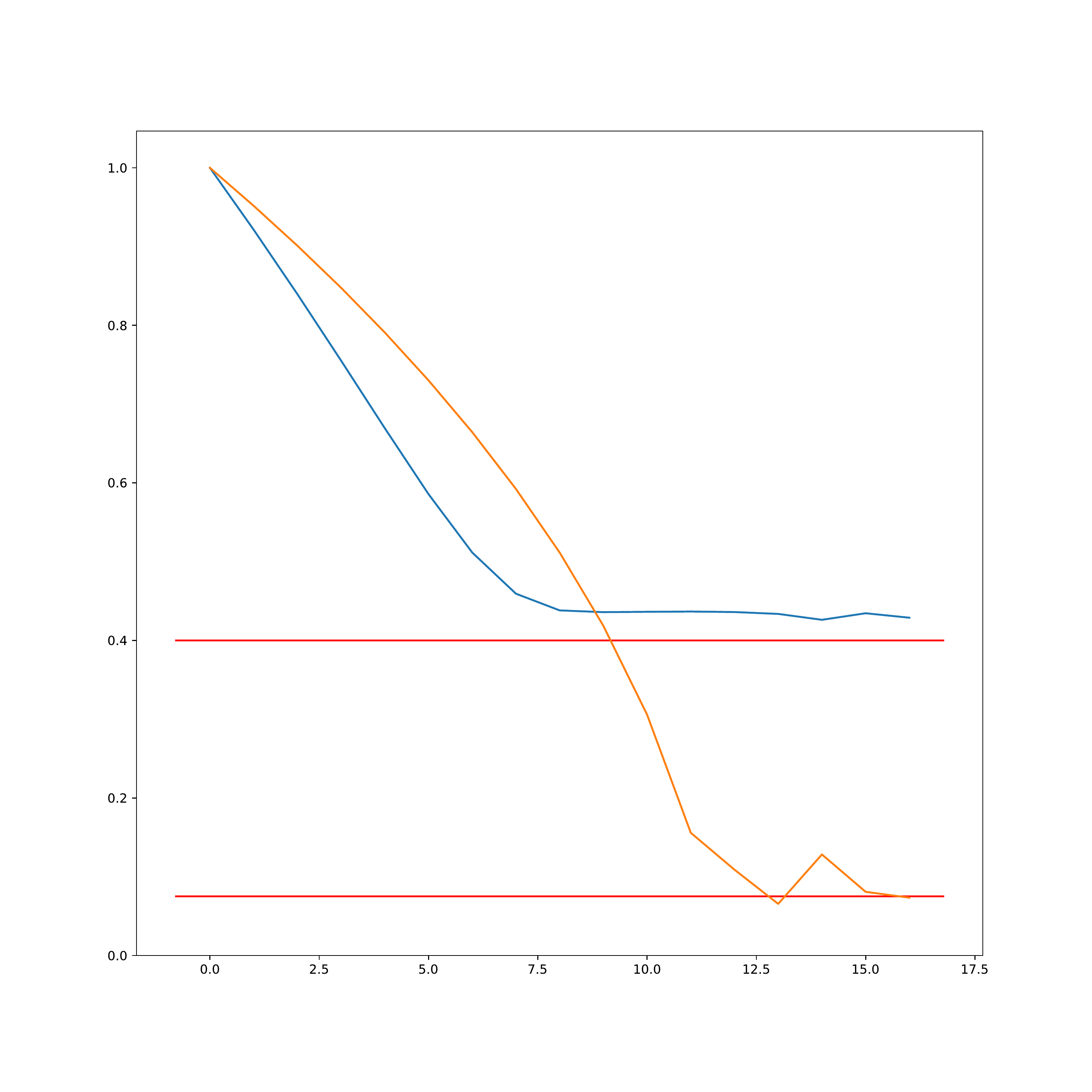}%
    \caption{(left) Estimated likelihood as a function of iterations of the iterative ML optimization. (right) Evolution of the variance and noise $\sigma$ during the iterations of the algorithm. Ground truth values in red.}
    \label{fig:likelihood}
\end{figure}

\subsection{Principal Components}
For a fixed point $v$ on the surfaces, Figure~\ref{fig:conditioned_samples} shows sample trajectories from the bridge process \eqref{eq:tildext} generated with a Hamiltonian MCMC sampler for approximate evaluation of the expectation in the density expression \eqref{eq:density_lim_E}. The corresponding samples from the latent process $x_t$ are shown in 
Figure~\ref{fig:conditioned_antidevelopment}. Notice how the endpoints of the latent process samples vary even though the trajectories on the surfaces always end at $v$. This effect is a direct consequence of non-zero curvature. Figure~\ref{fig:conditioned_antidevelopment_density} shows a density plot of the latent process, still conditioned on $v$. The mean latent path is plotted in blue in Figure~\ref{fig:conditioned_antidevelopment} and \ref{fig:conditioned_antidevelopment_density}, and the development of the mean path together with the parallel transported frame along the path are plotted on Figure~\ref{fig:conditioned_samples}. In Euclidean space, the mean path would be a straight line corresponding to the geodesics in Figure~\ref{fig:conditioned_samples}. While the mean latent paths for both surfaces clearly deviate from straight lines, the effect of the curvature is clearly more emphasized on the non-symmetric ellipsoid.

\subsection{Maximum Likelihood}
For the samples in Figure~\ref{fig:samples_mle}, we plot in Figure~\ref{fig:likelihood} the negative log-likelihood computed from a sample approximation of \eqref{eq:density_lim_E}, the estimated variance in the axis of major variation, and $\sigma$. The horizontal axis shows the evolution of the negative log-likelihood and estimated variance during the evolution of an iterative maximum likelihood optimization. The samples are generated with variance .4 in the major axis and $\sigma=.075$ noise. The algorithm makes repeated sample approximations of \eqref{eq:density_lim_E}, calculates the gradient and updates the parameters. As can be seen from the figure, the parameter estimates converges approximately to the true values.

\section{Conclusion and Outlook}
The probabilistic formulation in PPCA allows to generalize the PCA procedure to manifolds with a focus on data likelihoods in contrast to constructions of subspaces. This has previously been pursued with probabilistic PGA \cite{zhang_probabilistic_2013}. Here, we provide a generalization based on a different probability model using stochastic flows in the frame bundle and related fiber bundles. The main feature of the model is the intrinsic definition that does not refer to a linear tangent space approximation, is infinitesimal in modelling stochastic differential flows, and focuses on the generated likelihood and density instead of squared Riemannian distances. The model uses fiber bundle geometry that reveals important geometric information about the construction. As an example, the non-integrability of the horizontal subbundle is directly related to the curvature of the manifold. Instead of truncating the non-closure of the bracket of the horizontal basis fields to provide a submanifold, the construction allows the diffusion to spread into higher-dimensional subspaces. The data manifold is thereby not linearized and the curvature preserved in the analysis of the data.

The construction is based on the anisotropic normal distributions defined in \cite{sommer_anisotropic_2015,sommer_modelling_2017}. In addition to the presented PCA formulation, a regression model based on these distributions has been presented in \cite{kuhnel_stochastic_2017}. We hope in future work to be able to use this and similar geometric constructions that preserve the nonlinear nature of the data space to generalize more statistical procedures to analysis of manifold valued data in intrinsic ways.

\section*{Acknowledgments}
The work was supported by the Danish Council for Independent Research,
and the CSGB Centre for Stochastic Geometry and Advanced Bioimaging funded by a grant from the Villum foundation. The research was
partially performed at the Mathematisches Forschungsinstitut Oberwolfach (MFO), 2014 and 2018.

\bibliographystyle{spmpsci}      

\bibliography{ss}

\begin{thebibliography}{10}
\providecommand{\url}[1]{{#1}}
\providecommand{\urlprefix}{URL }
\expandafter\ifx\csname urlstyle\endcsname\relax
  \providecommand{\doi}[1]{DOI~\discretionary{}{}{}#1}\else
  \providecommand{\doi}{DOI~\discretionary{}{}{}\begingroup
  \urlstyle{rm}\Url}\fi

\bibitem{arnaudon_geometric_2018}
Arnaudon, A., Holm, D.D., Sommer, S.: A {Geometric} {Framework} for
  {Stochastic} {Shape} {Analysis}.
\newblock accepted for Foundations of Computational Mathematics,
  arXiv:1703.09971 [cs, math]  (2018)

\bibitem{delyon_simulation_2006}
Delyon, B., Hu, Y.: Simulation of conditioned diffusion and application to
  parameter estimation.
\newblock Stochastic Processes and their Applications \textbf{116}(11),
  1660--1675 (2006).
\newblock \doi{10.1016/j.spa.2006.04.004}

\bibitem{eltzner_torus_2015}
Eltzner, B., Huckemann, S., Mardia, K.V.: Torus {Principal} {Component}
  {Analysis} with an {Application} to {RNA} {Structures}.
\newblock arXiv:1511.04993 [q-bio, stat]  (2015).
\newblock \urlprefix\url{http://arxiv.org/abs/1511.04993}.
\newblock ArXiv: 1511.04993

\bibitem{elworthy_geometric_1988}
Elworthy, D.: Geometric aspects of diffusions on manifolds.
\newblock In: P.L. Hennequin (ed.) École d'Été de {Probabilités} de
  {Saint}-{Flour} {XV}–{XVII}, 1985–87, no. 1362 in Lecture {Notes} in
  {Mathematics}, pp. 277--425. Springer Berlin Heidelberg (1988).
\newblock \urlprefix\url{http://link.springer.com/chapter/10.1007/BFb0086183}

\bibitem{fletcher_principal_2004-1}
Fletcher, P., Lu, C., Pizer, S., Joshi, S.: Principal geodesic analysis for the
  study of nonlinear statistics of shape.
\newblock Medical Imaging, IEEE Transactions on  (2004).
\newblock \doi{10.1109/TMI.2004.831793}

\bibitem{frechet_les_1948}
Frechet, M.: Les éléments aléatoires de nature quelconque dans un espace
  distancie.
\newblock Ann. Inst. H. Poincaré \textbf{10}, 215--310 (1948)

\bibitem{hsu_stochastic_2002}
Hsu, E.P.: Stochastic {Analysis} on {Manifolds}.
\newblock American Mathematical Soc. (2002)

\bibitem{huckemann_intrinsic_2010}
Huckemann, S., Hotz, T., Munk, A.: Intrinsic shape analysis: {Geodesic} {PCA}
  for {Riemannian} manifolds modulo isometric {Lie} group actions.
\newblock Statistica Sinica \textbf{20}(1), 1--100 (2010)

\bibitem{jung_analysis_2012}
Jung, S., Dryden, I.L., Marron, J.S.: Analysis of principal nested spheres.
\newblock Biometrika \textbf{99}(3), 551--568 (2012).
\newblock \doi{10.1093/biomet/ass022}

\bibitem{kolar_natural_1993}
Kolář, I., Slovák, J., Michor, P.W.: Natural {Operations} in {Differential}
  {Geometry}.
\newblock Springer Berlin Heidelberg, Berlin, Heidelberg (1993).
\newblock \urlprefix\url{http://link.springer.com/10.1007/978-3-662-02950-3}

\bibitem{kuhnel_differential_2017}
Kühnel, L., Arnaudon, A., Sommer, S.: Differential geometry and stochastic
  dynamics with deep learning numerics.
\newblock arXiv:1712.08364 [cs, stat]  (2017).
\newblock \urlprefix\url{http://arxiv.org/abs/1712.08364}.
\newblock ArXiv: 1712.08364

\bibitem{kuhnel_stochastic_2017}
Kühnel, L., Sommer, S.: Stochastic {Development} {Regression} on {Non}-linear
  {Manifolds}.
\newblock In: Information {Processing} in {Medical} {Imaging}, Lecture {Notes}
  in {Computer} {Science}, pp. 53--64. Springer, Cham (2017).
\newblock \doi{10.1007/978-3-319-59050-9_5}

\bibitem{marchand_conditioning_2011}
Marchand, J.L.: Conditioning diffusions with respect to partial observations.
\newblock arXiv:1105.1608 [math]  (2011).
\newblock \urlprefix\url{http://arxiv.org/abs/1105.1608}.
\newblock ArXiv: 1105.1608

\bibitem{mok_differential_1978}
Mok, K.P.: On the differential geometry of frame bundles of {Riemannian}
  manifolds.
\newblock Journal Fur Die Reine Und Angewandte Mathematik \textbf{1978}(302),
  16--31 (1978).
\newblock \doi{10.1515/crll.1978.302.16}

\bibitem{pennec_barycentric_2016}
Pennec, X.: Barycentric {Subspace} {Analysis} on {Manifolds}.
\newblock arXiv:1607.02833 [math, stat]  (2016).
\newblock \urlprefix\url{http://arxiv.org/abs/1607.02833}.
\newblock ArXiv: 1607.02833

\bibitem{roweis_em_1998}
Roweis, S.: {EM} {Algorithms} for {PCA} and {SPCA}.
\newblock In: Proceedings of the 1997 {Conference} on {Advances} in {Neural}
  {Information} {Processing} {Systems} 10, {NIPS} '97, pp. 626--632. MIT Press,
  Cambridge, MA, USA (1998)

\bibitem{sommer_horizontal_2013}
Sommer, S.: Horizontal {Dimensionality} {Reduction} and {Iterated} {Frame}
  {Bundle} {Development}.
\newblock In: Geometric {Science} of {Information}, {LNCS}, pp. 76--83.
  Springer (2013)

\bibitem{sommer_diffusion_2014}
Sommer, S.: Diffusion {Processes} and {PCA} on {Manifolds}.
\newblock Mathematisches Forschungsinstitut Oberwolfach
  https://www.mfo.de/document/1440a/OWR\_2014{\textbackslash}\_44.pdf (2014).
\newblock \urlprefix\url{https://www.mfo.de/document/1440a/OWR_2014\_44.pdf}

\bibitem{sommer_anisotropic_2015}
Sommer, S.: Anisotropic {Distributions} on {Manifolds}: {Template} {Estimation}
  and {Most} {Probable} {Paths}.
\newblock In: Information {Processing} in {Medical} {Imaging}, \emph{Lecture
  {Notes} in {Computer} {Science}}, vol. 9123, pp. 193--204. Springer (2015)

\bibitem{sommer_evolution_2015}
Sommer, S.: Evolution {Equations} with {Anisotropic} {Distributions} and
  {Diffusion} {PCA}.
\newblock In: F.~Nielsen, F.~Barbaresco (eds.) Geometric {Science} of
  {Information}, no. 9389 in Lecture {Notes} in {Computer} {Science}, pp.
  3--11. Springer International Publishing (2015).
\newblock \doi{10.1007/978-3-319-25040-3_1}

\bibitem{sommer_anisotropically_2016}
Sommer, S.: Anisotropically {Weighted} and {Nonholonomically} {Constrained}
  {Evolutions} on {Manifolds}.
\newblock Entropy \textbf{18}(12), 425 (2016).
\newblock \doi{10.3390/e18120425}

\bibitem{sommer_diffusion_2018}
Sommer, S.: Diffusion {Bridge} {Simulation} on {Nonlinear} {Manifolds}.
\newblock in preparation  (2018)

\bibitem{sommer_bridge_2017}
Sommer, S., Arnaudon, A., Kuhnel, L., Joshi, S.: Bridge {Simulation} and
  {Metric} {Estimation} on {Landmark} {Manifolds}.
\newblock In: Graphs in {Biomedical} {Image} {Analysis}, {Computational}
  {Anatomy} and {Imaging} {Genetics}, Lecture {Notes} in {Computer} {Science},
  pp. 79--91. Springer (2017).
\newblock \doi{10.1007/978-3-319-67675-3_8}

\bibitem{sommer_modelling_2017}
Sommer, S., Svane, A.M.: Modelling anisotropic covariance using stochastic
  development and sub-{Riemannian} frame bundle geometry.
\newblock Journal of Geometric Mechanics \textbf{9}(3), 391--410 (2017).
\newblock \doi{10.3934/jgm.2017015}

\bibitem{the_theano_development_team_theano:_2016}
Team, T.T.D.: Theano: {A} {Python} framework for fast computation of
  mathematical expressions.
\newblock arXiv:1605.02688 [cs]  (2016).
\newblock ArXiv: 1605.02688

\bibitem{tipping_probabilistic_1999}
Tipping, M.E., Bishop, C.M.: Probabilistic {Principal} {Component} {Analysis}.
\newblock Journal of the Royal Statistical Society. Series B \textbf{61}(3),
  611--622 (1999)

\bibitem{zhang_probabilistic_2013}
Zhang, M., Fletcher, P.: Probabilistic {Principal} {Geodesic} {Analysis}.
\newblock In: {NIPS}, pp. 1178--1186 (2013)

\end{thebibliography}


\end{document}